%% file: paperII.tex
\documentclass[11pt,final]{./Siam/siamltex}

\usepackage{hyperref}

\usepackage{pifont}
\newcommand{\xmark}{\ding{55}}
\usepackage{float}

\usepackage{enumerate}

\usepackage{amsmath,tikz}

\usetikzlibrary{arrows.meta}
\usepackage{xcolor}

\input{commands}

\newtheorem{example}[theorem]{Example}
\newtheorem{assumption}[theorem]{Assumption}

\newcommand{\rmdiv}{\mathrm{div}}

\newcommand{\Gammain}{{\Gamma_{in}}}
\newcommand{\Ltwoin}{L^2_{V_n}(\Gammain)}
\newcommand{\ltwoin}[1]{\norm{#1}_{\Ltwoin}}

\newcommand{\Gammaout}{{\Gamma_{out}}}
\newcommand{\Ltwoout}{L^2_{V_n}(\Gammaout)}
\newcommand{\ltwoout}[1]{\norm{#1}_{\Ltwoout}}

\newcommand{\delTh}{\partial \calT_h}

\newcommand{\delTho}{\partial \mathring{\calT}_h}
\newcommand{\LtwoDelTho}{L^2_{V_n}(\delTho)}
\newcommand{\ltwoDelTho}[1]{\norm{#1}_{\LtwoDelTho}}

\usepackage{enumitem}
\setlist[itemize]{leftmargin=*}
\setlist[enumerate]{leftmargin=*}
\begin{document}

\title{Asymptotic Robustness of HDG Schemes
  for the Convection-Diffusion Equation}

\author{
  Lucas Bouck\thanks{Center for Computational Biology, 
  Flatiron Institute, New York, NY 10010.} \and
  Noel J. Walkington\thanks{Department of Mathematical Sciences,
    Carnegie Mellon University, Pittsburgh, PA 15213;
    supported in part by NSF grant RTG-2342349, and
    grant DMS-2424556 while the author was in residence at
    ICERM in Providence, RI.}}

\date{\today}

\maketitle

\begin{abstract}
  The convection-diffusion problem has been studied extensively as an
  important model for convection-dominated problems. The
  problem becomes more challenging when the convecting velocity
  is no longer Lipschitz and merely has Sobolev regularity.  In this
  paper, we use Boyer's extension of the DiPerna--Lions theory for
  the transport equation to study the vanishing diffusion limit of
  discretizations of the convection-diffusion equation with Sobolev
  velocity. In the case of flux-adapted boundary conditions, we prove
  that a continuous Galerkin discretization converges strongly to the
  renormalized solution of the transport equation as the diffusion
  coefficient and mesh size go to zero, without requiring any
  relationship between them. We refer to this unconditional strong
  convergence as asymptotic robustness.  In the case of Dirichlet
  boundary conditions, we use the renormalization theory and barrier
  arguments to prove that a piecewise constant hybridizable
  discontinuous Galerkin scheme is asymptotically robust, while the
  continuous Galerkin scheme may fail to converge.  Our results
  indicate that the asymptotic robustness of schemes for the
  convection-diffusion problem depends on the realization of the
  boundary conditions as opposed to internal layers.
\end{abstract}


% \begin{MSCcodes}
% Codes here
% \end{MSCcodes}

\section{Introduction}
Continuous Galerkin (CG) and discontinuous Galerkin (DG)
approximations of solutions to the transport and
convection-diffusion equations
\begin{align}
  u + \bfV.\nabla u &= f,
  && \text{ in } \Omega,  \label{eqn:strong-transport} \\
  u_\epsilon + \bfV.\nabla u_\epsilon - \epsilon \Delta u_\epsilon & = f,
  && \text{ in } \Omega,
\label{eqn:strong-advection-diffusion}
\end{align}
are analyzed. Here $\Omega\subset \mathbb{R}^d$ is a bounded domain
and $\epsilon > 0$.

Below we show that if the boundary conditions are ``adapted'' (see
equation~\eqref{eqn:adapted-bc} below), then most finite element
schemes for the convection-diffusion problem converge in $\Ltwo$ to a
solution of the transport equation independently of how the mesh and
diffusion parameters $(h,\epsilon)$ converge to $(0,0)$.  We will
refer to this as ``asymptotic robustness'' with respect to the
indicated parameters.  It is well known that this is not the case if
Dirichlet boundary conditions are prescribed for the
convection-diffusion problem; indeed, in the absence of restrictions
relating $h$ and $\epsilon$, many schemes exhibit oscillations and so,
at best, converge weakly in $\Ltwo$.  Below we show that piecewise
constant hybridizable discontinuous Galerkin (HDG) approximations of
the convection-diffusion equation with Dirichlet boundary data are
asymptotically robust; that is, they converge strongly in $\Ltwo$ to a
solution of the transport equation independently of how
$(h,\epsilon) \rightarrow (0, 0)$. Similar notions of
asymptotic robustness appear as asymptotic
compatibility for nonlocal equations~\cite{tian2014asymptotically},
robustness in linear elasticity~\cite{babuvska1992locking}, and
asymptotic-preserving schemes for hyperbolic
equations~\cite{jin1999efficient}.

When the velocity field is not Lipschitz, solutions of the
equations~\eqref{eqn:strong-transport}
and~\eqref{eqn:strong-advection-diffusion} may develop layers (or
singularities) in $\Omega$ in addition to the classical boundary
layers. One unexpected outcome of the results below is that the
boundary layers appear to be the main obstacle to strong convergence
of (stable) numerical schemes. Specifically, if the boundary
conditions for the convection-diffusion equation are ``adapted'' in
the sense that
\begin{equation} \label{eqn:adapted-bc}
  u_\epsilon \bfV.\bfn +\epsilon \dbydp{u_\epsilon}{n}
  = u_0 \bfV.\bfn
  \quad \text{ on } \Gamma_{in},
  \qquad \text{ and } \qquad
   \dbydp{u_\epsilon}{n} =0
  \quad \text{ on }
  \partial\Omega\setminus \Gamma_{in},
\end{equation}
Theorem~\ref{thm:adapted-bc-conv} below shows that the classical CG finite element
scheme is asymptotically robust; moreover, the proof of strong
convergence carries over {\em mutatis mutandis} to most stable
numerical schemes. Here $\Gamma_{in}$ is the inflow boundary,
$$
\Gamma_{in} = \{x \in \partial \Omega \sst \bfV.\bfn < 0\}
\qquad \text{ and } \qquad
\dbydp{u_\epsilon}{n} = (\nabla u_\epsilon).\bfn,
$$
where $\bfn$ denotes the unit outward-pointing normal to $\Omega$.
Example~\ref{ex:cg-counter-example} below shows that CG approximations
of the convection-diffusion equation with Dirichlet boundary data,
\begin{equation}\label{eqn:dirichlet-bc}
u|_{\partial \Omega} = u_0,
\end{equation}
 may not converge strongly; in addition,
their weak limits may be solutions of the transport equation with
boundary data, $u|_{\Gamma_{in}} \neq u_0$.

The arguments below carry over to the corresponding evolution problems
when $\bfV$ is independent of time.
Equations~\eqref{eqn:strong-transport}
and~\eqref{eqn:strong-advection-diffusion} arise when the evolution
problems
\begin{align}
  \partial_t u + \bfV.\nabla u
  &= f \quad \text{ in } (0,T)\times \Omega, \label{eqn:evol-transport} \\
  \partial_t u_\epsilon + \bfV.\nabla u_\epsilon - \epsilon \Delta u_\epsilon
  &= f \quad \text{ in } (0,T)\times \Omega, \label{eqn:evol-advdiff}
\end{align}
are approximated using the implicit Euler time-stepping scheme;
$\partial_t u(t) \approx (u(t+\tau) - u(t))/\tau$.  If
$\bfV \in \Hone^d$ with $\linf{\rmdiv(\bfV)} < \infty$, semigroup
theory shows that discrete solutions of the time-dependent problem
will converge in $C[0,T; \Ltwo]$ as
$(\tau, h, \epsilon) \rightarrow (0,0,0)$ if numerical solutions of
the corresponding resolvent equations, \eqref{eqn:strong-transport}
and~\eqref{eqn:strong-advection-diffusion}, converge strongly in
$\Ltwo$ as
$(h,\epsilon) \rightarrow (0,0)$~\cite{brezis1972convergence,
  crandall1971generation, itokappel1998trotter}.  In this context the
spatial terms are scaled by the time step $\tau$, so it suffices to
consider the resolvent equations with velocity fields satisfying
$\linf{\rmdiv(\bfV)} < 2$.

\subsection{Renormalization, Vanishing Diffusion Limit, and Main Contributions}

Our work utilizes Boyer's extension of the
DiPerna--Lions~\cite{diperna1989ordinary} theory to accommodate
nonhomogeneous boundary data~\cite{boyer2005trace}.
DiPerna--Lions introduced the notion of renormalized solution
to the transport equation, which provides a notion of well-posedness
for the transport equation when $\bfV$ is not Lipschitz and one cannot
apply the method of characteristics. If $\bfV$ is smooth and $u$ is a
smooth solution to $\partial_t u + \bfV.\nabla u = f$ on
$(0,T)\times \mathbb{R}^d$, then
$\partial_t \beta(u) + \bfV.\nabla \beta(u) = f\beta'(u)$ for
$\beta\in C^1(\mathbb{R})$. DiPerna and Lions showed that this holds
true if $u \in \Lp$ is a distributional
solution of the transport equation with $\bfV \in \Wonepo$
for all $\beta\in C^1(\mathbb{R})$ with suitable growth.

Formulas for the $\Lp$ norms of distributional solutions of the
transport equation can be developed by selecting $\beta(s) = |s|^p$
and exploiting the renormalization property. Strong convergence of
weakly converging (approximate) solutions
$\{u_\epsilon\}_{\epsilon > 0}$ can frequently be established by
showing convergence of their norms,
$\lp{u_\epsilon} \rightarrow \lp{u}$.  This line of argument is used
below to establish strong convergence of numerical solutions
computed with CG/DG Galerkin schemes and a piecewise constant DG/HDG
scheme. This leads to the following question.  \vspace{.25cm}
\begin{center}
  \parbox{.95\textwidth}{ \textbf{Question}: Let
    $\{u_{h\epsilon}\}_{h,\epsilon > 0}$ be numerical approximations of
    the convection-diffusion
    problem~\eqref{eqn:strong-advection-diffusion} with mesh size
    $h$. Let $u\in \Ltwo$ be the weak solution
    to~\eqref{eqn:strong-transport} with inflow boundary
    specified. Does
    $$
    \lim_{(h,\epsilon)\to(0,0)} u_{h\epsilon} = u
    $$
    strongly in $\Ltwo$ independently of how $(h,\epsilon) \to (0,0)$?
  }
\end{center}
\vspace{.25cm}

While it is often simple to prove that numerical approximations will
converge weakly to a distributional solution of the transport equation
as $(h,\epsilon)\to(0,0)$, preserving the boundary conditions and
strong $\Ltwo$ convergence are more difficult problems. The following
example illustrates that the answer to the above question is negative
for numerical approximations computed using a CG scheme with Dirichlet
boundary conditions.

\begin{example}\label{ex:cg-counter-example}
Consider the boundary value problem on $\Omega = (0,1)$ with Dirichlet boundary conditions:
\begin{equation}\label{eqn:1d-bvp}
\begin{aligned}
u' - \epsilon u'' &= 0, \text{ in } (0,1), \quad u(0) = 0,\quad u(1) =1.
\end{aligned}
\end{equation}
The exact solution to~\eqref{eqn:1d-bvp} is
$u_\epsilon(x) = \frac{e^{x/\epsilon} - 1}{e^{1/\epsilon} - 1}.$
As $\epsilon\to0$, the strong $\Ltwo$ limit of $u_\epsilon$ is
$u(x) = 0$, which solves the transport equation
\begin{equation}\label{eqn:1d-transport}
\begin{aligned}
u' &= 0, \text{ in } (0,1), \quad u(0) = 0.
\end{aligned}
\end{equation}
Let $u_{h\epsilon}$ be the CG piecewise linear solution
of~\eqref{eqn:1d-bvp} on a uniform grid with spacing $h = 1/N$.  If
$h,\epsilon\to0$ with $h^2/2\epsilon \to a > 0$, there are two
subsequences of $\{u_{h\epsilon}\}_{h,\epsilon > 0}$ (not relabeled)
that have the following weak $\Ltwo$ limits:
$$
u_{h\epsilon} \rightharpoonup u^\pm = \frac{1}{1 \pm e^{2/a}}.
$$
Moreover, $u^\pm$ are the only cluster points of the sequence
$\{u_{h\epsilon}\}$, and there is no cluster point that is the correct
transport solution.  While the weak limits of $u_{h\epsilon}$ are
distributional solutions to the transport
equation~\eqref{eqn:1d-transport}, they are nonzero constants so do
not satisfy the inflow boundary condition $u(0) = 0$. This
oscillatory failure of the standard CG discretization is
an instance of the general mechanism for convection-dominated
discretizations described by
Bacuta~\cite{bacuta2026explaining}. Figure~\ref{fig:counter-example}
displays the discrete CG solution with this oscillatory behavior, the
true transport solution, and the incorrect weak limit.

\begin{figure}
\begin{center}
\includegraphics[width=.6\textwidth]{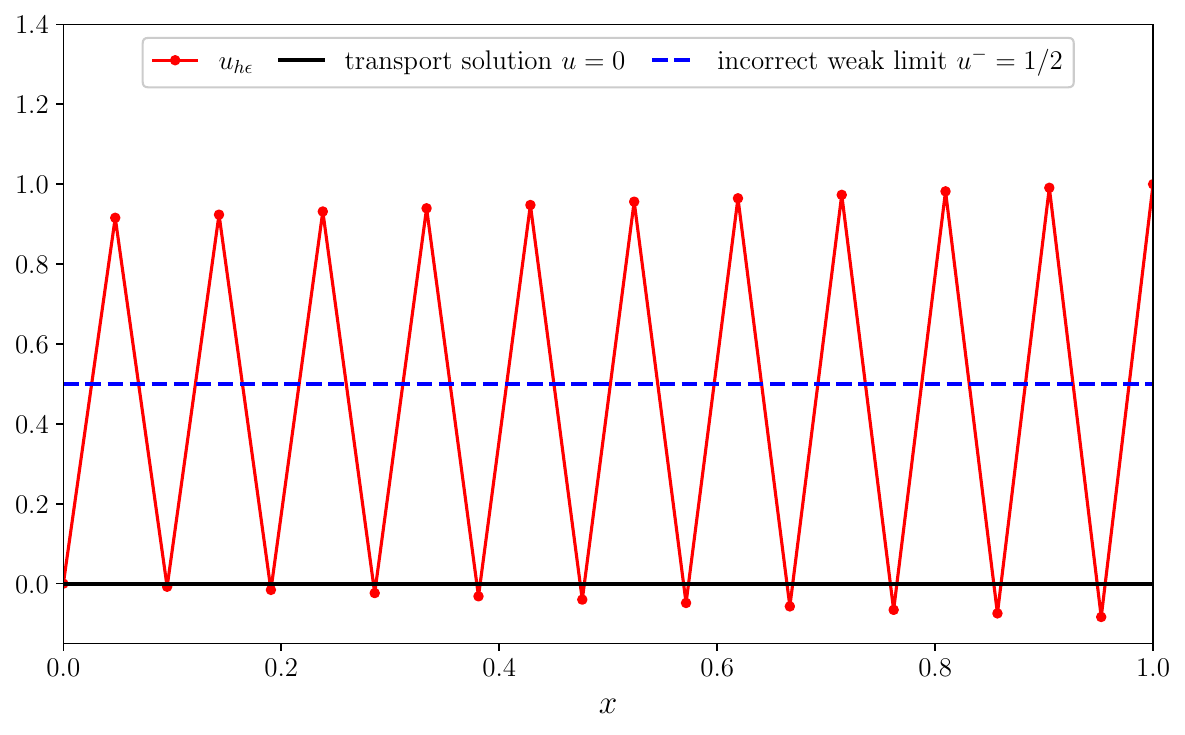}
\caption{Discrete solution $u_{h\epsilon}$ with $\epsilon = 10^{-4}$
  and $N = 21$.  If $h^2/\epsilon\to\infty$, then for odd $N$ the
  nodal values converge to $0$ on even nodes and $1$ on odd nodes, so
  $u_{h\epsilon} \rightharpoonup \half$, which does not agree with the
  correct limit $u(x) = 0$.}
  \label{fig:counter-example}
\end{center}
\end{figure}
\end{example}

Our results show that the asymptotic robustness of schemes depends
greatly on the boundary conditions. Boyer's trace
theorems~\cite{boyer2005trace} for functions having convected
derivatives in $\Lp$ are used below to identify the limits of
numerical solutions. We show that the standard CG scheme converges
strongly in $\Ltwo$ when the boundary conditions are adapted. Although
CG fails with Dirichlet boundary conditions, we show this can be
mitigated for a piecewise constant DG/HDG scheme in
Section~\ref{sec:dirichlet-bc}. Table~\ref{tab:summary} summarizes our
answers to this question for different boundary conditions and
schemes.
\begin{table}[H]
\begin{center}
\begin{tabular}{c|c|c}
\parbox{.25\textwidth}{Method/BC}&   \parbox{.3\textwidth}{Adapted Flux BC~\eqref{eqn:adapted-bc}} &  \parbox{.25\textwidth}{Dirichlet BC~\eqref{eqn:dirichlet-bc}}  \\
 \hline
 \parbox{.25\textwidth}{Continuous Problem}   &   \parbox{.33\textwidth}{\checkmark\; (see arguments in~\cite{boyer2005trace})} &  \parbox{.33\textwidth}{\checkmark\; Thm.~\ref{thm:dirichlet-conv-cont}}   \\
 \parbox{.25\textwidth}{Continuous Galerkin} & \parbox{.33\textwidth}{\checkmark\; Thm.~\ref{thm:adapted-bc-conv}} &  \parbox{.33\textwidth}{\xmark\; when $h^2 \gg \epsilon$ \, Ex.~\ref{ex:cg-counter-example} }  \\
 \parbox{.33\textwidth}{Piecewise Constant DG/HDG (Sec.~\ref{sec:pw-const-scheme})} & \parbox{.33\textwidth}{\checkmark\; (see proof of Thm.~\ref{thm:adapted-bc-conv}) } & \parbox{.33\textwidth}{\checkmark\; Thm.~\ref{thm:dirichlet-conv-hdg} }
\end{tabular}
\end{center}
\caption{Summary of known results and our results. Does
  $u_{h\epsilon}\to u$ strongly in $L^2(\Omega)$ in the dual limit
  $h,\epsilon\to0$ for any relationship between $h$ and $\epsilon$?}
\label{tab:summary}
\end{table}

\subsection{Related Literature}
Numerics for the transport and convection-diffusion equations have
been extensively studied and the literature is vast. Since an
exhaustive review of this is not feasible, we focus on the works
particularly close to ours. For a broader view of DG methods and
convection-dominated problems, we point the reader to the
books~\cite{di2011mathematical, riviere2008discontinuous,
  roos2008robust} and the references therein.

Relatively few works consider approximation of low-regularity
renormalized solutions to the evolution transport
problem~\eqref{eqn:evol-transport} with Sobolev or $\mathrm{BV}$
velocity field $\bfV$. Walkington~\cite{walkington2005convergence}
proved strong $L^2$ convergence in space and time for DG
approximations assuming the velocity
$\bfV \in L^1[0,T; H^1_0(\Omega)^d]$, and Boyer
\cite{boyer2012analysis} proved strong convergence of a finite volume
scheme assuming $\bfV \in L^1[0,T; W^{1,1}(\Omega)^d]$ using his
results in~\cite{boyer2005trace}. Later, Aguillon and
Boyer~\cite{aguillon2018error} prove $\mathcal{O}(\sqrt{h})$ error
estimates in $L^1$ between an upwind solution and the solution of the
transport equation~\eqref{eqn:evol-transport} assuming
$\bfV\in \mathrm{BV}[0,T; L^{\infty}(\Omega)^d]\cap L^\infty[0,T;
W^{1,\infty}(\Omega)^d]$, and Schlichting and
Seis~\cite{schlichting2018analysis} prove $\mathcal{O}(\sqrt{h})$
error estimates in a weaker metric assuming
$\bfV\in L^1[0,T;W^{1,p}(\Omega)^d]$.

For stationary problems such as~\eqref{eqn:strong-transport}
and~\eqref{eqn:strong-advection-diffusion}, Ern and
Guermond~\cite{ern2006discontinuous} developed a general stability and
convergence theory of DG methods for Friedrichs' systems. Part of
their theory involved defining trace operators similar to those of
Boyer (see~\cite[Section~2.1.3]{di2011mathematical} for additional
discussion).  A consequence of their general theory is that a DG
solution for the transport equation will converge strongly in $\Ltwo$
as $h\to0$ to the solution of~\eqref{eqn:strong-transport} provided
$\bfV$ is sufficiently regular and the inflow and outflow boundaries
are well separated. Another consequence of their theory is error
estimates and convergence as $h\to0$ of a DG method for the
convection-diffusion problem~\eqref{eqn:strong-advection-diffusion}
for fixed $\epsilon>0$. However, they do not study the dual limit
$(h,\epsilon)\to(0,0)$. Bause and Knabner~\cite{bause2002uniform}
studied a Lagrange--Galerkin scheme for the evolution
convection-diffusion problem~\eqref{eqn:evol-advdiff}, and proved
error estimates between $u_\epsilon$ and $u_{h\epsilon}$ in $\Ltwo$
that are independent of $\epsilon$. However, their analysis requires
sufficient smoothness of $\bfV$ to apply the method of
characteristics.

Hybridized formulations for convection-diffusion problems have also
been studied extensively. Egger and Sch\"oberl~\cite{egger2010hybrid}
propose a hybrid mixed DG method that combines a mixed discretization
of the diffusion with an upwind DG discretization of the transport
part, with the two components coupled through hybridization. By
construction, the method reduces to the corresponding mixed method in
the purely diffusive limit and to the corresponding DG method in the
hyperbolic limit. Assuming $\bfV$ is a piecewise constant field, Egger
and Sch\"oberl prove~\cite[Theorem~3.14]{egger2010hybrid} energy-norm
a priori error estimates between the discrete solution,
$u_{h\epsilon}$, and the convection-diffusion solution, which imply
$ \Vert u_{h\epsilon} - u_\epsilon\Vert_{\Ltwo} \leq C
h^s(\sqrt{\epsilon} +\sqrt{h
  |\bfV|})|u_\epsilon|_{H^{1+s}(\mathcal{T}_h)}, $ where $s>1/2$ and
$H^{1+s}(\mathcal{T}_h)$ should be understood in a broken sense. In
later work, Fu, Qiu, and Zhang~\cite{fu2015analysis} build
on~\cite{egger2010hybrid} and design an equivalent HDG scheme by
carefully choosing the stabilization. For
$\bfV\in W^{1,\infty}(\Omega)^d$ and their HDG method with piecewise
linear functions, they prove a priori error estimates on general
meshes of the form
$ \Vert u_{h\epsilon} - u_\epsilon\Vert_{\Ltwo}\leq C h(\epsilon
h^{-1/2} + \sqrt{\epsilon} +\sqrt{h})| u_\epsilon |_{H^2(\Omega)} $ in
the convection-dominated regime. They also prove error estimates with
optimal rates in $h$ for special meshes following techniques developed
in~\cite{cockburn2008optimal}.

A priori error estimates involving Sobolev norms of the solution do
not immediately imply that $u_{h\epsilon}$ converges to the transport
solution $u$ strongly in $\Ltwo$ for any sequence
$(h,\epsilon)\to(0,0)$ because the energy estimates implicitly rely on
the regularity of the convection-diffusion solution $u_\epsilon$.
With internal or boundary layers, the norms of the
convection-diffusion solution scale like
$|u_\epsilon|_{H^{1+s}(\Omega)} \approx \epsilon^{-s - 1/2}$ for
$0\leq s\leq 1$.  As a result, showing $u_{h\epsilon}\to u$ in $\Ltwo$
using the a priori estimates from~\cite{egger2010hybrid,
  fu2015analysis} requires $h/\epsilon \to 0$.  To the best of our
knowledge, existing works do not establish, or disprove, in the case
of $\bfV\in \Hone^d$, strong convergence in $\Ltwo$ to the
renormalized inflow solution without requiring a relationship between
$h$ and $\epsilon$, a property we will refer to as asymptotic
robustness.  Our results show, moreover, that this convergence depends
decisively on the realization of the boundary conditions.

\section{Preliminaries}
We consider a bounded domain $\Omega\subset \mathbb{R}^d$ with
polygonal boundary, and denote the inner product on $\Ltwo$ by
$(u,v)$.  The velocity fields are assumed to satisfy
$\bfV\in H^1(\Omega)^d$ with $\rmdiv(\bfV) \in \Linf$. The
corresponding inflow and outflow boundaries are
$$
\Gammain := \{x\in \partial\Omega : \bfV.\bfn < 0\},
\qquad \text{ and } \qquad
\Gammaout := \{x\in \partial\Omega : \bfV.\bfn > 0\}.
$$
The boundary may then be partitioned as
$\partial \Omega = \Gammain \cup \Gamma_0 \cup \Gammaout$, where
$\Gamma_0$ is the characteristic boundary. The weak formulations will involve
boundary integrals with weight $|\bfV.\bfn|$ on the inflow and outflow
boundaries.  The weighted $L^2$ spaces are denoted by
$$
\Ltwoin = \{u:\Gammain \rightarrow \Re
\sst \int_{\Gammain} \!\! u^2 \, |\bfV.\bfn| \, dS < \infty\},
$$
and similarly for $\Ltwoout$,
with the corresponding inner products denoted by
$(u,v |\bfV.\bfn|)_{\Gammain}$ and $(u,v |\bfV.\bfn|)_{\Gammaout}$.

Convergence of numerical approximations will be established for a
regular sequence of triangulations of $\Omega$ denoted by
$\{\mathcal{T}_{h}\}_{h > 0}$. For each mesh, its $(d-1)$-skeleton is
denoted by $\partial \mathcal{T}_h$, and its interior $(d-1)$-skeleton
by $\partial \mathring{\mathcal{T}}_h$. The space of discontinuous
piecewise polynomials of degree $k$ on $\mathcal{T}_h$ is denoted by
$\mathcal{P}_k(\mathcal{T}_h)$, and the space of discontinuous
piecewise polynomials of degree $k$ defined on the skeleton is denoted
by $\mathcal{P}_k(\partial \mathcal{T}_h)$. The $L^2$-projection
operators $P_h: L^2(\Omega)\to \mathcal{P}_0(\mathcal{T}_h)$ and
$\Phat_h: L^2(\Omega)\to \mathcal{P}_0(\partial \mathcal{T}_h)$ are
$$
P_h v|_K = \frac{1}{|K|} \int_K v, \quad K \in \calT_h,
\quad \text{ and } \quad
\Phat_h v|_k = \frac{1}{|k|} \int_k v, \quad k \in \partial \calT_h,
$$
and $\mathbf{P}_h: L^2(\Omega)^d\to \mathcal{P}_k( \calT_h)^d$ is the
vector-valued analog.  The Lagrange interpolant that maps
$u\in C(\Omegabar)$ onto $\mathcal{P}_k(\mathcal{T}_h)\cap C(\Omegabar)$
is denoted by $I_h$.

\subsection{Boyer's Extension of DiPerna--Lions Theory}
The seminal paper of DiPerna and Lions~\cite{diperna1989ordinary}
established existence and uniqueness of solutions of the transport
equation\footnote{DiPerna--Lions and Boyer consider the evolution
  equation~\eqref{eqn:evol-transport}; their results extend directly to the
  stationary problem.} with velocity
fields satisfying $\bfV \in \Wonepo^d$ with
$\linf{\rmdiv(\bfV)} < \infty$. A key property of the solutions was
that they could be ``renormalized''; that is, if $u, f \in \Ltwo$ and
$$
u + \bfV.\nabla u = f \quad \text{ in } \calD'(\Omega)
\qquad \text{ then } \qquad
\beta'(u) u + \bfV.\nabla \beta(u) = \beta'(u) f, 
\quad \text{ in } \calD'(\Omega),
$$
for all $\beta \in C^1(\Re)$ with suitable growth.

When $\bfV \in \Wonepo^d$, no boundary conditions are required for the
transport equation. In a series of papers, Boyer extended the
DiPerna--Lions theory to admit more general velocity fields, and
developed the corresponding trace theory for functions $u \in \Lp$
with distributional derivative
$\bfV.\nabla u \in \Lp$~\cite{boyer2005trace}. For smooth functions
the trace of a function, $\gamma(u) = u|_{\partial\Omega}$, is
characterized by
\begin{equation} \label{eqn:trace}
\left(\phi \bfV. \bfn,\gamma(u) \vph\right)_{\partial \Omega}
= (\bfV.\nabla u , \phi)  + (u, \mathrm{div}(\phi \bfV)),
\qquad \phi \in C^\infty(\Omegabar).
\end{equation}
The following result of Boyer shows that this characterization of the
trace is well defined when $\bfV$ and $u$ are less regular, and will
be used to identify boundary conditions satisfied by limits of the
numerical approximations to the transport and convection-diffusion
equations.

\begin{proposition} [Boyer~\cite{boyer2005trace,boyer2012mathematical}]
  \label{prop:trace}
  Let $\Omega \subset \Re^d$ be a Lipschitz domain and
  $\bfV \in \Hone^d$ with $\linf{\rmdiv(\bfV)} < \infty$. If
  $u\in L^2(\Omega)$ has distributional derivative
  $\bfV.\nabla u \in L^2(\Omega)$, then equation~\eqref{eqn:trace}
  characterizes a measurable function, $\gamma(u)$, on
  $\partial \Omega$ which is uniquely defined almost everywhere with
  respect to the measure $|\bfV. \bfn| \, dS$; in addition, if
  $\gamma(u)\in \Ltwoin$, then $\gamma(u)\in \Ltwoout$.
\end{proposition}

In general, $\gamma(u) \not\in L^2_{V_n}(\partial\Omega)$~\cite[Remark~2.6]{di2011mathematical}; an essential property of this
trace is that if $\gamma(u) \in \Ltwoin$, then this regularity gets
``pushed/transported'' to the outflow.  This proposition is the
technical tool required to extend the DiPerna--Lions theory and
prescribe boundary data for solutions of the transport equation when
$\bfV \in \Hone^d$.

\begin{theorem}[Well-posedness and renormalization property]
  \label{thm:boyer-transport}
  Let $\bfV \in \Hone^d$ with $\linf{\rmdiv(\bfV)} < 2$,
  $f \in \Ltwo$, and $u_0 \in \Ltwoin$. Then there exists a unique
  pair $(u, \ut) \in \Ltwo \times \Ltwoout$ such that 
  $$
  (u, v) - ( u, \rmdiv(v \bfV)) + (\ut, v \bfV.\bfn)_\Gammaout
  = (f,v) - (u_0, v \bfV.\bfn)_\Gammain,
  \quad v\in C^1(\Omegabar);
  $$
  in particular, $\gamma(u) = \ut$ a.e.\ on $\Gammaout$ and $\gamma(u) = u_0$
  on $\Gammain$. In addition, $(u,\gamma(u))$ is a renormalized solution
  in the following sense. If $\beta \in C^1(\mathbb{R})$ satisfies
  $|\beta(s)| \leq C(1+|s|^{\min(2, 1 + 2/d)})$ and
  $|\beta'(s)| \leq C(1+|s|)$, then
  $$
  (u\beta'(u), v) - ( \beta(u), \rmdiv(v \bfV))
  + (\beta(\tilde{u}), v \bfV.\bfn)_\Gammaout
  = (f\beta'(u),v)
  - (\beta(u_0), v \bfV.\bfn)_\Gammain,
  $$
  for all $v\in C^1(\Omegabar)$.
\end{theorem}

The restriction on the growth of $\beta$ guarantees that the term
$(\beta(u), \rmdiv(v \bfV))$ is integrable. Under additional
hypotheses, the restrictions on the growth of $\beta$ and the regularity of
the test functions $v$ can be relaxed as in the following corollary.

\begin{corollary}\label{cor:test-by-u}
  Let $\bfV, u_0, f$ satisfy the hypotheses of Theorem~\ref{thm:boyer-transport}, and
  $(u,\gamma(u))\in \Ltwo \times \Ltwoout$ be the solution of
  the corresponding transport equation. Then
  \begin{equation}\label{eqn:test-by-u}
    \int_\Omega  \big(1- \frac{1}{2} \rmdiv( \bfV)\big) u^2 dx
    + \frac{1}{2} \int_\Gammaout \gamma(u)^2 |\bfV.\bfn| dS
    = (f, u)
    - \half (u_0^2, \bfV.\bfn)_\Gammain.
  \end{equation}
  If, in addition, $u \in \Lq$ with $q > 4d/(d+1)$ and
  $\phi \in \Woneinf$, then
  \begin{equation}\label{eqn:test-by-u-phi}
    \int_\Omega  \big(\phi - \frac{1}{2} \rmdiv(\phi \bfV)\big) u^2 \, dx 
    + \half (\gamma(u)^2, \phi \bfV.\bfn)_\Gammaout.
    = (f, \phi u) - \half (u_0^2, \phi \bfV.\bfn)_\Gammain.
  \end{equation}
\end{corollary}

\begin{proof}
  Approximate $\beta(s) = s^2/2$ by
  $$
  \beta_n(s) = \left\{
    \begin{array}{cc}
      s^2 / 2 & |s| \leq n, \\ n |s| - n^2/2 & n \leq |s|;
    \end{array}
  \right.
  \quad \text{ in particular } \quad
  0 \leq \beta_n(s) \leq s^2 / 2
  \text{ and }
  |\beta'_n(s)| \leq |s|.
  $$
  Then $\beta_n \in C^1(\mathbb{R})$ has linear growth and
  bounded derivative, so if $v \in C^1(\Omegabar)$, the renormalization
  property shows
  \begin{equation*}
    (u\beta'_n(u), v) - ( \beta_n(u), \rmdiv(v \bfV))
    + (\beta_n(\gamma(u)), v \bfV.\bfn)_\Gammaout
    = (f \beta'_n(u),v) - (\beta_n(u_0), v \bfV.\bfn)_\Gammain.
  \end{equation*}
  When $v = 1$, the term $(\beta_n(u), \rmdiv(\bfV))$ is dominated by
  $\linf{\rmdiv(\bfV)} u^2/2$, and $\beta_n(\gamma(u))$ and
  $\beta'_n(u)$ converge pointwise to, and are dominated by,
  $\gamma(u)^2/2$ and $u$ respectively. We may then pass to the limit by
  the Lebesgue Dominated Convergence Theorem to get equation~\eqref{eqn:test-by-u}.

  To establish equation~\eqref{eqn:test-by-u-phi}, note that
  $\rmdiv(v \bfV) \in \Linf$ when $u \in \Lq$, so it serves as a
  dominating function. We may then pass to the limit in the above to
  get
  \begin{equation} \label{eqn:usqv}
    \int_\Omega \left(v  - \half \rmdiv(v \bfV) \right) u^2
    + \half \left(\gamma(u)^2, v \bfV.\bfn \vph\right)_\Gammaout
    = (f u,v) - \half \left(u_0^2, v \bfV.\bfn \vph\right)_\Gammain.
  \end{equation}
  If $\phi \in W^{1,\infty}(\Omega)$, let
  $\phi_\eta \in W^{1,\infty}(\Re^d)$ denote the mollification of the
  Lipschitz extension of $\phi$ to $\Re^d$ with parameter $\eta > 0$. Then
  $$
  \linf{\phi_\eta} \leq \linf{\phi},
  \qquad \text{ and } \qquad
  \linf{\nabla \phi_\eta} \leq \linf{\nabla \phi},
  $$
  so $\phi_\eta(x) \rightarrow \phi(x)$ uniformly on $\Omegabar$,
  and $\nabla\phi_\eta(x) \rightarrow \nabla \phi(x)$ for a.e.\
  $x \in \Omega$; in particular,
  $\rmdiv(\phi_\eta \bfV) \rightarrow \rmdiv(\phi \bfV)$
  at a.e.\ $x \in \Omega$.
  Since $\rmdiv(\phi_\eta \bfV) u^2$ is dominated by
  $C (|\bfV| + |\rmdiv(\bfV)|) u^2$, equation~\eqref{eqn:test-by-u-phi} follows upon setting
  $v = \phi_\eta$ in equation~\eqref{eqn:usqv} and passing to the
  limit using the Dominated Convergence
  Theorem.  
\end{proof}

If $u_n \weak u$ in $\Ltwo$, it can be difficult to determine if
$\gamma(u_n)$ converges to $\gamma(u)$. The following squeeze lemma will
be used to identify the limiting trace in this context.

\begin{lemma} \label{lem:squeezeTrace}
  Let $\Omega \subset \Re^d$ be a Lipschitz domain and
  $\bfV \in \Hone^d$ with $\linf{\rmdiv(\bfV)} < \infty$, and let
  $u, v_\pm \in \Ltwo$ have distributional derivatives
  $\bfV.\nabla u, \bfV.\nabla v_\pm \in \Ltwo$.  If
  $\Gamma \subset \partial \Omega$ and $\omega \subset \Omegabar$ is
  an open neighborhood of $\Gamma$ and
  $$
  v_- \leq u \leq v_+
  \quad \text{ on } \omega,
  \quad \text{ then } \quad
  \gamma(v_-) \leq \gamma(u) \leq \gamma(v_+)
  \quad \text{ on } \Gamma.
  $$
\end{lemma}

If $v \in \Hone$, it follows directly from the definition that
$\gamma(v) = v$ is the usual Sobolev trace, so if $v_+ = v_- = v$,
it is immediate that $\gamma(u) = v$. A key observation is that
sets of the form $\{u \in L^2(\omega) \sst v_- \leq u \leq v_+\}$
are weakly closed, so weak limits of functions in such a set have
traces satisfying $\gamma(v_-) \leq \gamma(u) \leq \gamma(v_+)$.

\begin{proof}[Sketch of Proof]
This result follows from Boyer's construction of the trace, which proceeds
in the following three steps~\cite{boyer2012mathematical}.
\begin{itemize}
\item A partition of unity $\{\phi_i\}_{i=1}^N$ is used to write $u$
  as the sum of functions $\phi_i u$, each of which has nonzero
  boundary values on a portion of $\partial \Omega$ which can be
  represented as the graph of a Lipschitz function.

\item The functions are mollified with a kernel with shifted origin.

\item The trace is obtained as the limit of the mollified functions.
\end{itemize}
The lemma then follows since each of these steps preserves both
linearity and ordering, and at each point, $x \in \partial \Omega$, the
trace depends only upon the function values in a neighborhood of $x$.
\end{proof}

\section{Adapted Boundary Conditions}
This section considers classical finite element approximations of
the convection-diffusion equation,
$$
u_\epsilon - \epsilon \Delta u_\epsilon + \bfV.\nabla u_\epsilon = f
\qquad \text{ in } \Omega,
$$
with boundary conditions
$$
(u_\epsilon \bfV + \epsilon \nabla u_\epsilon).\bfn = u_0 \bfV.\bfn
\quad \text{ on } \Gammain,
\qquad \text{ and } \qquad
\bfn.\nabla u_\epsilon =0  \quad \text{ on } \Gammaout.
$$
The following assumptions set the context for the numerical
approximation of this problem.

\begin{assumption}\label{as:adapted-bc}
Assume $\Omega \subset \Re^d$ is a bounded Lipschitz domain and:
\begin{enumerate}
\item $\bfV \in \Hone^d$ with $\linf{\rmdiv(\bfV)} < 2$.

\item $f\in \Ltwo$ and $u_0\in \Ltwoin$.

\item $\{\calT_h\}_{h > 0}$ is a regular family of triangulations of
  $\Omega$, parameterized by
  $h = \max_{K \in \calT_h} \mathrm{diam}(K)$, that also triangulates
  $\Gammain$ and $\Gammaout$.
\end{enumerate}
\end{assumption}

Classical finite element methodology approximates the weak statement: find
$u_\epsilon \in \Hone$ such that
$$
(u_\epsilon,v) + \epsilon (\nabla u_\epsilon, \nabla v)
- (u_\epsilon, \rmdiv(v \bfV)) + (u_\epsilon, v \bfV.\bfn)_\Gammaout
= (f,v) - (u_0, v \bfV.\bfn)_\Gammain,
\quad v \in \Hone.
$$
Setting the test function $v = u_\epsilon$ and integrating by parts
yields the classical a priori bound:
\begin{multline} \label{eqn:uepsBound}
  \ltwo{u_\epsilon}^2
  - \half \left(u_\epsilon^2, \rmdiv(\bfV) \vph\right)
  + \epsilon \ltwo{\nabla u_\epsilon}^2 \\
  + \half \ltwoin{u_\epsilon-u_0}^2
  + \half \ltwoout{u_\epsilon}^2 
  = (f,u_\epsilon) + \half \ltwoin{u_0}^2.
\end{multline}

The following theorem establishes strong convergence of Galerkin
approximations of the convection-diffusion problem using the classical
finite element spaces (CG schemes). The proof extends to essentially
all finite element schemes that exhibit an energy estimate of this
form.

\begin{theorem}\label{thm:adapted-bc-conv}
  Let Assumption~\ref{as:adapted-bc} hold, and set
  $U_h = \{u_h \in \Hone \sst u_h|_K \in \calP_k(K), \,\, K \in
  \calT_h\}$ with integer $k \geq 1$ fixed. For each $\epsilon > 0$ let
  $u_{h\epsilon} \in U_h$ satisfy
  \begin{equation} \label{eqn:uhepsWeak}
    (u_{h\epsilon},v_h) + \epsilon (\nabla u_{h\epsilon}, \nabla v_h)
    - (u_{h\epsilon}, \rmdiv(v_h \bfV))
    + (u_{h\epsilon}, v_h \bfV.\bfn)_\Gammaout
    = (f,v_h) - (u_0, v_h \bfV.\bfn)_\Gammain,
  \end{equation}
  for all $v_h \in U_h$.  Then
  $\lim_{(h,\epsilon) \rightarrow (0,0)} u_{h\epsilon} = u$ in
  $\Ltwo$, where $u$ is the (renormalized) solution of the transport
  equation
  $$
  u + \bfV.\nabla u = f, \quad \text{ in } \calD'(\Omega),
  \quad \text{ with } \quad \gamma(u)|_\Gammain = u_0.
  $$
  In addition, $\sqrt{\epsilon} \ltwo{\nabla u_{h\epsilon}} \rightarrow 0$
  and
  $\norm{u_{h\epsilon} - \gamma(u)}_{L^2_{V_n}(\partial \Omega)} \rightarrow 0$.
  
\end{theorem}

\begin{proof}
  The proof naturally breaks down into the following four steps.

{\em1. Weak Limit:} As for the continuous problem, setting
    $v_h = u_{h\epsilon}$ in equation~\eqref{eqn:uhepsWeak} gives the
    bound
    $$
    \ltwo{u_{h\epsilon}} + \sqrt{\epsilon} \ltwo{\nabla u_{h\epsilon}}
    + \ltwoin{u_0 - u_{h\epsilon}} + \ltwoout{u_{h\epsilon}}
    \leq C \left(\ltwo{f} + \ltwoin{u_0} \vph\right).
    $$
    Upon passing to a subsequence (not relabeled), it follows that
    there exists
    $(u,\ut) \in \Ltwo \times L^2_{V_n}(\partial \Omega)$
    for which $u_{h\epsilon} \weak u$ in $\Ltwo$ and
    $u_{h\epsilon}|_{\partial \Omega} \weak \ut$ in
    $L^2_{V_n}(\partial \Omega)$. 

    {\em2. Limit Equation:} Fix $v \in \calD(\Omega)$ and select
    $v_h \in U_h$ in~\eqref{eqn:uhepsWeak} to be the Lagrange
    interpolant of $v$ to get
    $$
    (u_{h\epsilon},v_h) + \epsilon (\nabla u_{h\epsilon}, \nabla v_h)
    - (u_{h\epsilon}, \rmdiv(v_h \bfV))
    = (f,v_h) .
    $$
    Since $v_h \rightarrow v$ in $\Hone$ as $h \rightarrow 0$,
    $\sqrt{\epsilon}\nabla u_{h\epsilon}$ is uniformly bounded in
    $\Ltwo^d$, and $u_{h\epsilon} \weak u$ in $\Ltwo$, it follows that
    $u \in \Ltwo$ satisfies $u + \bfV.\nabla u = f$ in
    $\calD'(\Omega)$. In particular, the distributional derivative
    $\bfV.\nabla u$ belongs to $\Ltwo$, and
    Proposition~\ref{prop:trace} shows that $u$ has a trace,
    $\gamma(u)$.

{\em3. Boundary Values:} Let $v \in C^1(\Omegabar)$ and select
    $v_h \in U_h$ in~\eqref{eqn:uhepsWeak} to be the Lagrange interpolant
    of $v$ to get
    $$
    (u_{h\epsilon},v_h) + \epsilon (\nabla u_{h\epsilon}, \nabla v_h)
    - (u_{h\epsilon}, \rmdiv(v_h \bfV))
    + (u_{h\epsilon}, v_h \bfV.\bfn)_\Gammaout
    = (f,v_h) - (u_0, v_h \bfV.\bfn)_\Gammain.
    $$
    Since $\sqrt{\epsilon}\nabla u_{h\epsilon}$ is uniformly bounded
    in $\Ltwo^d$, and $v_h\to v$ strongly in $\Hone$ and in the uniform
    norm, we may pass to the limit to derive
    $$
    (u,v) - (u, \rmdiv(v \bfV)) + (\ut, v \bfV.\bfn)_\Gammaout
    = (f,v) - (u_0, v \bfV.\bfn)_\Gammain,
    \quad v \in C^1(\Omegabar).
    $$
    Rearranging the above equation and using the definition of
    $\gamma(u)$ yields
    $$
    (\ut, v \bfV.\bfn)_\Gammaout + (u_0, v \bfV.\bfn)_\Gammain
    = (f,v) -(u,v) + (u, \rmdiv(v \bfV)) 
    = (\gamma(u),v \bfV.\bfn)_{\partial \Omega}.
    $$
    It follows that $\gamma(u)|_\Gammain = u_0 \in \Ltwoin$ and
    $\gamma(u)|_\Gammaout = \ut \in \Ltwoout$. In particular, since
    the solution of the transport equation is unique, the sequence
    $\{u_{h\epsilon}\}_{h,\epsilon > 0}$ converges weakly without
    passage to a subsequence. Moreover, since $u$ is a renormalized
    solution and $u_0, f$ satisfy the hypotheses of Theorem~\ref{thm:boyer-transport},
    it follows from Corollary~\ref{cor:test-by-u} that
    $$
    \ltwo{u}^2 - \half\left(u^2, \rmdiv(\bfV)\right)
    +  \half\ltwoout{\gamma(u)}^2
    = (f,u) +\half  \ltwoin{u_0}^2.
    $$
    
    {\em4. Strong Convergence:} Since
    $u_{h\epsilon}\rightharpoonup u$ in $\Ltwo$, strong convergence
    follows from convergence of the $\Ltwo$ norm. Recall the energy
    estimate~\eqref{eqn:uepsBound}, and write the same estimate for
    the CG solution:
    \begin{align*}
      \ltwo{u_{h\epsilon}}^2
      - \half \left(u_{h\epsilon}^2, \rmdiv(\bfV) \vph\right)
      &+ \epsilon \ltwo{\nabla u_{h\epsilon}}^2 \\
      &+ \half \ltwoin{u_{h\epsilon}-u_0}^2
        + \half \ltwoout{u_{h\epsilon}}^2
        = (f,u_{h\epsilon}) + \half \ltwoin{u_0}^2.
    \end{align*}
    Omitting the positive terms
    $\epsilon \ltwo{\nabla u_{h\epsilon}}^2 + \half
    \ltwoin{u_{h\epsilon}-u_0}^2$ on the left-hand side, we calculate
    \begin{align*}
      \limsup_{(h,\epsilon) \rightarrow (0,0)}
      \Big\{
      \ltwo{u_{h\epsilon}}^2 - \half (u_{h\epsilon}^2, \rmdiv(\bfV))
      &+ \half \ltwoout{u_{h\epsilon}}^2  \Big\} \\
      &\leq \limsup_{(h,\epsilon) \rightarrow (0,0)}
          (f,u_{h\epsilon}) +\half \ltwoin{u_0}^2 \\
      &= (f,u) + \half\ltwoin{u_0}^2 \\
      &= \ltwo{u}^2
          - \half \left(u^2, \rmdiv(\bfV) \vph\right)
          + \half \ltwoout{\gamma(u)}^2 .
    \end{align*}
    Since $1 - \half \rmdiv(\bfV)$ is bounded away from zero, it
    follows that $u_{h\epsilon} \to u$ in $\Ltwo$, and
    $u_{h\epsilon}|_\Gammaout \to \gamma(u)$ in $\Ltwoout$. A squeeze
    argument with the omitted terms then shows
    $\sqrt{\epsilon}\nabla u_{h\epsilon} \to0$ in $\Ltwo^d$ and
    $u_{h\epsilon}\to u_0$ in $\Ltwoin$.
\end{proof}
\subsection{Computational Example}\label{sec:comp-ex-adapted}
\begin{figure}[H]
\begin{center}
  \includegraphics[height=5.5cm]{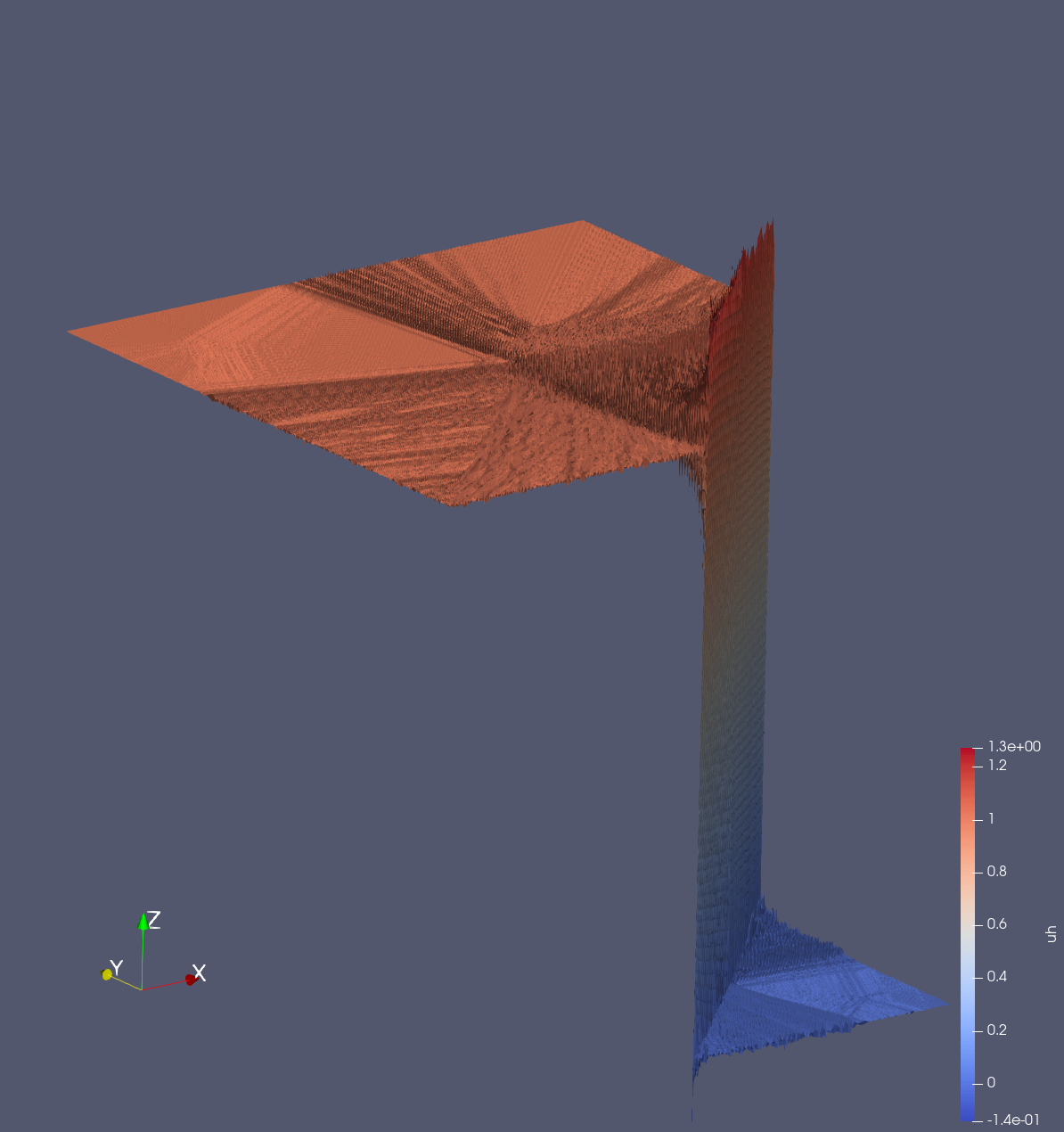}
  \includegraphics[height=5.5cm]{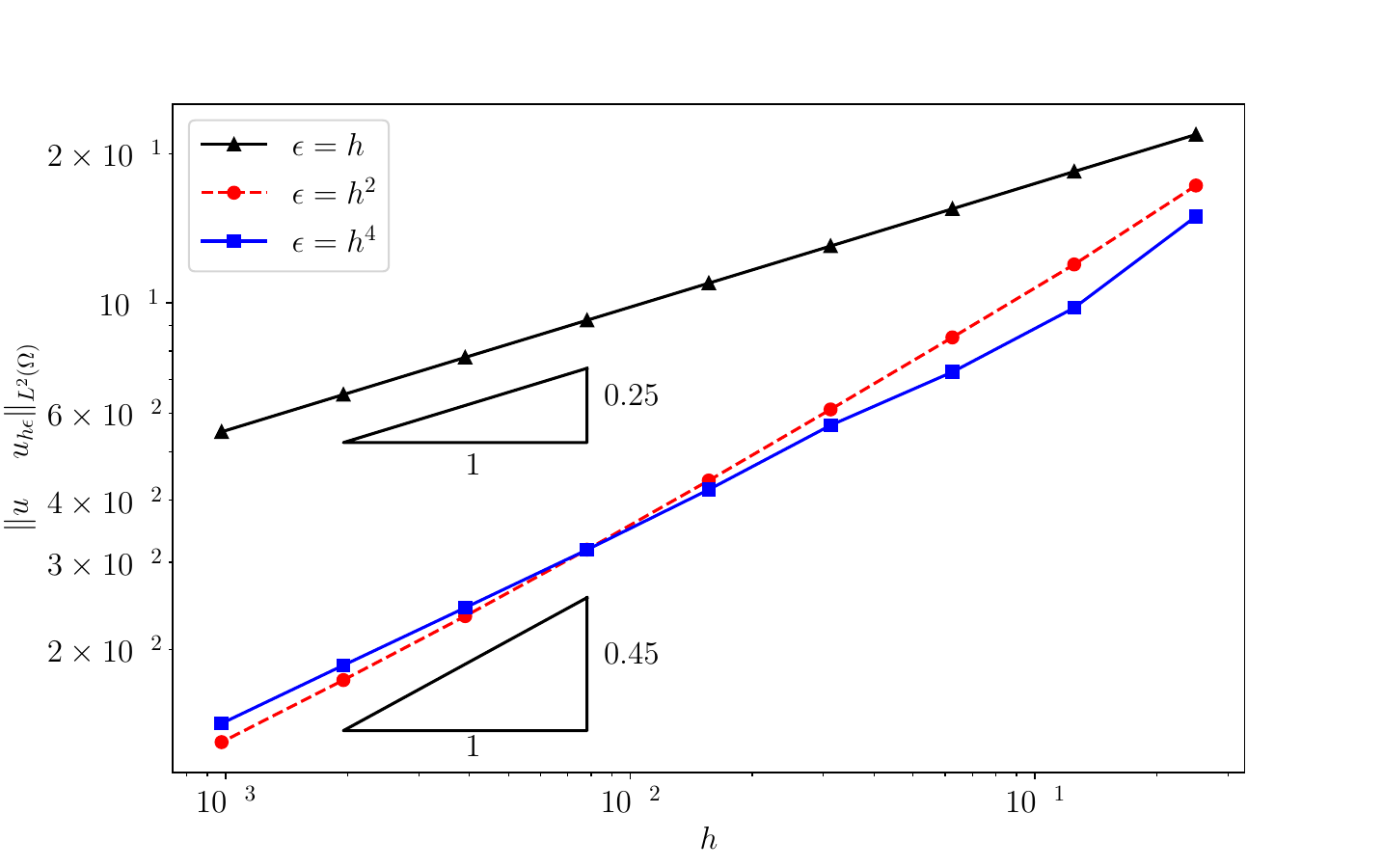}
  \caption{(Left) CG solution $u_{h\epsilon}$ with $h=2^{-10}$,
    $\epsilon = h^2$. (Right) $\Vert u - u_{h\epsilon}\Vert_{\Ltwo}$ vs.\
    $h$ for $\epsilon =h, h^2, h^4$.
  }
\label{fig:cg-convergence}
\end{center}
\end{figure}

To illustrate the results in Theorem~\ref{thm:adapted-bc-conv}, we
implemented the following example in NGSolve~\cite{schoberl2014c++} with domain, velocity, and manufactured
solution:
$$
\Omega =(0,1)^2, \quad \bfV
= (1,1)^T, \quad u(x) =
\begin{cases} 1, &  x_1-x_2 > \frac{1}{2}\\
  0, & \text{ otherwise}
\end{cases}.
$$ 
The right-hand side $f$ can then be constructed by evaluating
$u+\bfV.\nabla u$. We take $h = 2^{-\ell}$ for $\ell = 2, \ldots, 10$
and test the CG scheme with $\epsilon = h, h^2, h^4$. For all three
relationships,
$\Vert u - u_{h\epsilon}\Vert_{\Ltwo}$ decreases with consistent
rates. Figure~\ref{fig:cg-convergence} shows the CG solution
$u_{h\epsilon}$ and a log-log plot of
$\Vert u - u_{h\epsilon}\Vert_{\Ltwo}$. The CG solution oscillates
at the level of the mesh, but the $\Ltwo$ norm of these oscillations vanishes in the limit
as $h,\epsilon\to 0$ with an empirical rate of convergence of
$O(\sqrt{h} +\epsilon^{1/4})$.

\section{Dirichlet Boundary Conditions}\label{sec:dirichlet-bc}
This section considers the convection-diffusion equation with Dirichlet
boundary data,
\begin{equation} \label{eqn:strong0}
u_\epsilon - \epsilon \Delta u_\epsilon + \bfV.\nabla u_\epsilon = f
\qquad \text{ in } \Omega,
\qquad \text{ with } \quad
u_\epsilon|_{\partial\Omega} = 0.
\end{equation}
For simplicity, homogeneous Dirichlet data is assumed on the entire
boundary. The analysis below also holds if Neumann or Robin
boundary data is specified on $\Gammaout$, and if the Dirichlet data on $\Gammain$
is the trace of a function in $\Hone$.

Solutions of~\eqref{eqn:strong0}, $u_\epsilon \in \Honeo$, satisfy
\begin{equation} \label{eqn:weak0}
(u_\epsilon,v) + \epsilon (\nabla u_\epsilon, \nabla v)
- (u_\epsilon, \rmdiv(v \bfV))
= (f,v), \qquad v \in \Honeo.
\end{equation}
As for adapted boundary conditions, setting the test function
$v = u_\epsilon$ gives the classical energy estimate
\begin{equation} \label{eqn:energy0}
  \ltwo{u_\epsilon}^2
  - \half \left(u_\epsilon^2, \rmdiv(\bfV) \vph\right)
  + \frac{\epsilon}{2} \ltwo{\nabla u_\epsilon}^2
  = (f, u_\epsilon).
\end{equation}
It follows that $\{u_\epsilon\}_{\epsilon > 0}$ is bounded in $\Ltwo$,
and passing to the limit in~\eqref{eqn:weak0} shows that weak limits
satisfy the transport equation $u + \bfV.\nabla u = f$ in
$\calD'(\Omega)$. The major difficulty is to verify $\gamma(u) = 0$
on $\Gammain$. The maximum principle is used to establish this.

\subsection{Continuous Problem}
This subsection identifies weak limits of solutions
$\{u_\epsilon\}_{\epsilon > 0}$ of~\eqref{eqn:strong0} using a barrier argument.

\begin{definition}\label{def:barrier}
  A {\em barrier} at $\Gammain \subset \partial \Omega$ of equation~\eqref{eqn:strong0} is a function $\psi \in Lip(\omega)$ defined on a
  Lipschitz subset $\omega = \Omega\cap O$, where $O$ is an open
  neighborhood of $\Gammain$, satisfying
  \begin{itemize}
  \item $\psi(x) \geq 0$ on $\omega$, and $\psi(x) = 0$ if $x \in \Gammain$.

  \item $\psi - \epsilon \Delta \psi + \bfV.\nabla \psi \geq |f|$ on
    $\omega$ in the sense of distributions.

  \item There exists $\delta > 0$ such that $\psi(x) \geq \delta$ for
    $x \in \partial\omega \setminus \partial \Omega$.
  \end{itemize}
\end{definition}

Frequently the distance function to $\Gammain$ is a barrier function,
and this is the case we will consider below. Alternative constructions
of sub- and supersolutions of the transport equation can be found in~\cite[Appendix~A]{ayuso2009discontinuous} and~\cite[Sec.~2]{devinatz1974asymptotic}, which use the method of
characteristics.

\begin{example}
  Let $\Omega\subset \mathbb{R}^2$ be a polygonal domain and let
  $\Omega$ be convex near $\Gammain$ with the maximum angle between $\Gammain$
  and $\partial\Omega \setminus\Gammain$ bounded by $\pi/2$, and
  $\bfV$ continuous in a neighborhood of $\Gammain$. For suitable
  right-hand sides $f$, the distance function
  $d(x) = \mathrm{dist}(x, \Gammain)$ is a barrier on
  $\omega = \Gammain \cup \{x \in \Omega \sst d(x) < \delta\}$.  A
  sketch of the argument is as follows:
  \begin{itemize}
  \item $d(x) \geq 0$ and vanishes on $\Gamma_{in}$.
    
  \item If $\Omega$ is convex near $\Gammain$, then $-\Delta d \geq 0$
    for $d(x) < \delta$ sufficiently small~\cite{armitage1985convexity}, and $\nabla d.\bfV = -\bfV.\bfn > 0$
    on $\Gammain$. If $\bfV$ is continuous it follows that
    $d - \epsilon \Delta d + \bfV.\nabla d \geq 0$ on a neighborhood
    of $\Gammain$. Upon scaling $d(.)$ by a positive constant, the
    left-hand side can frequently bound $|f|$.
    
  \item $d(x) = \delta$ on $\partial\omega \setminus \partial\Omega$.
  \end{itemize}
\end{example}

The convergence theorem is proved under the following assumptions on
the data and barrier function.

\begin{assumption}\label{as:dirichlet-bc}
  Assume $\Omega$ is a bounded Lipschitz domain and,
  \begin{itemize}
  \item $f \in \Linf$, $\bfV\in \Hone^d \cap \Linf^d$, and
    $\linf{\rmdiv(\bfV)} < 2$.

  \item There exists a barrier $\psi$ at $\Gammain$ of equation~\eqref{eqn:strong0} in the sense of Definition~\ref{def:barrier}.
  
\item There exists $\delta > 0$ such that the distance function
  $d_{out}(x) = \mathrm{dist}(x, \Gammaout)$ satisfies
  $\bfV.\nabla d_{out} \leq 0$ on
  $\{x \in \Omegabar \sst 0 \leq d_{out}(x) \leq \delta \}$.
  \end{itemize}
\end{assumption}

If $f\not\in \Linf$, a Lipschitz barrier function satisfying the
conditions in Definition~\ref{def:barrier} sometimes does not
exist.

\begin{example}
  If $f(x) = -(x+1) \log(x)$ and $\bfV(x) \equiv 1$, then
  $f \in L^2(0,1)$ and the solution to $u + V u' = f$ with $u(0)=0$ is
  $u(x) = -x \log(x) + 1 - e^{-x} \geq 0$. Since $u'(0) = \infty$
  there is no (Lipschitz) barrier.
\end{example}

The following elementary lemma is used in the proof of the
main theorem of this section.

\begin{lemma}\label{lem:strong-L2-conv}
  Let $\Omega \subset \Re^d$ be a Lipschitz domain,
  $\phi \in W^{1,\infty}(\Omega)$ satisfy $\phi> 0$ a.e.\ in
  $\Omega$, and let $u \in L^\infty(\Omega)$. If the sequence
  $\{u_h\}_{h>0} \subset \Ltwo$ satisfies 
  $\sup_{h>0} \Vert u_h\Vert_{L^\infty(\Omega)} < \infty$ and
  $\Vert \phi^{1/2} (u - u_h)\Vert_{\Ltwo} \to0$ as $h\to0$, then
  $u_h \to u$ in $\Ltwo$.
\end{lemma}

A proof of this lemma is provided in the Appendix.

\begin{theorem}\label{thm:dirichlet-conv-cont}
  Let Assumption~\ref{as:dirichlet-bc} hold, and set
  $\phi(x) = \min(\mathrm{dist}(x, \Gammaout), \delta)$, where
  $\delta$ is the constant in Assumption~\ref{as:dirichlet-bc}. If
  $u_\epsilon$ denotes the solution to~\eqref{eqn:weak0}, then
  $\lim_{\epsilon \rightarrow 0} \ltwo{u_\epsilon - u} = 0$ where $u$
  is the (renormalized) solution of the transport equation
  $$
  u + \bfV.\nabla u = f, \quad \text{ in } \calD'(\Omega),
  \quad \text{ with } \quad \gamma(u)|_\Gammain = 0.
  $$
  In addition,
  $\sqrt{\epsilon} \ltwo{\phi^{1/2} \nabla u_{\epsilon}} \rightarrow
  0$.
\end{theorem}

\begin{proof}
  The proof proceeds in steps similar to the proof of Theorem~\ref{thm:adapted-bc-conv}.

{\em1. Weak Limit:} The energy estimate~\eqref{eqn:energy0} shows
    $\{u_\epsilon\}_{\epsilon > 0}$ is bounded in $\Ltwo$ and, upon
    passage to a subsequence, there exists a weak limit
    $u_\epsilon \weak u$ in $\Ltwo$.
    
{\em2. Limit Equation:} Passing to the limit in the
    weak statement~\eqref{eqn:weak0}, with $v \in \Honeo$ fixed, shows
    $u + \bfV.\nabla u = f$ in $\calD'(\Omega)$.
    
    {\em3. Boundary Values:} The maximum principle for equation~\eqref{eqn:strong0} gives $ \linf{u_\epsilon} \leq \linf{f}$. It
    follows that the barrier can be scaled to get
    $\psi(x) \geq |u_\epsilon(x)|$ for
    $x \in \partial\omega \setminus \partial \Omega$. Then, in the
    sense of distributions,
    $$
    (u_\epsilon-\psi) - \epsilon \Delta (u_\epsilon-\psi)
    + \bfV.\nabla (u_\epsilon-\psi)
    = f - (\psi + \bfV.\nabla \psi - \epsilon \Delta \psi) \leq 0,
    \qquad \text{ on } \omega.
    $$
    Since $u_\epsilon - \psi \leq 0$ on $\partial \omega$, it follows
    from the maximum principle that $u_\epsilon \leq \psi$ on
    $\omega$. Combining this with the analogous argument for
    $u_\epsilon + \psi$ shows $-\psi \leq u_\epsilon \leq \psi$ on
    $\omega$. The set of functions
    $\{u \in L^2(\omega) \sst -\psi \leq u \leq \psi\}$, being closed
    and convex, is weakly closed. Consequently, weak limits of
    $\{u_\epsilon\}_{\epsilon > 0}$ satisfy
    $-\psi \leq u|_\omega \leq \psi$ on $\omega$. Lemma~\ref{lem:squeezeTrace}
    then shows $\gamma(u)|_\Gammain = \psi|_\Gammain = 0$.

    It follows that $\{u_\epsilon\}_{\epsilon > 0}$ converges weakly
    to the solution of $u + \bfV.\nabla u = f$ with trace
    $\gamma(u)|_\Gammain = 0$.  Since the solution of the transport
    equation is unique, the sequence $\{u_\epsilon\}_{\epsilon > 0}$
    converges weakly without passage to a subsequence.
    
    {\em4. Strong Convergence:} To establish strong convergence, let
    $\phi(x) = \min(\mathrm{dist}(x, \Gammaout), \delta)$.  Since $u$
    is a renormalized solution and the trace
    $\gamma(\phi u) \in L^2\left(\partial\Omega; |\bfV.\bfn|\right)$
    vanishes, it follows that
    $$
    \int_\Omega (\phi - (1/2) \rmdiv(\phi \bfV)) u^2 = \int_\Omega \phi f u,
    $$
    and setting $v = \phi u_\epsilon$ in the weak statement~\eqref{eqn:weak0} then gives
    $$
    \int_\Omega \left(\phi - (1/2) \rmdiv(\phi \bfV) \vph\right) u_\epsilon^2
    + \epsilon \nabla u_\epsilon . \nabla (\phi u_\epsilon)
    = \int_\Omega \phi f u_\epsilon.
    $$
    Then
    \begin{eqnarray*}
      \lefteqn{\limsup_{\epsilon \rightarrow 0}
      \Big\{
      \int_\Omega \left(\phi - (1/2) \rmdiv(\phi \bfV) \vph\right) u_\epsilon^2
      + \phi \epsilon |\nabla u_{\epsilon}|^2 
      + \epsilon u_\epsilon \nabla u_\epsilon.\nabla \phi \Big\} }
      \hspace*{2.5in} \\
      &=& \limsup_{\epsilon \rightarrow 0}
      (\phi f,u_{\epsilon}) \\
      &=& (\phi f,u) \\
      &=& \int_\Omega \left(\phi - (1/2) \rmdiv(\phi \bfV) \vph\right) u^2.
    \end{eqnarray*}
    The cutoff function $\phi$ was constructed to satisfy
    $$
    \phi(x) > 0, \qquad
    \bfV.\nabla \phi \leq 0, \qquad 
    |\nabla \phi| \leq 1,
    \quad \text{ on } \Omega, \quad \text{ and } \quad
    \phi(x) = 0 \quad \text{ on } \Gammaout, 
    $$
    so that there is a constant $c>0$ such that
    $$
    \phi - (1/2) \rmdiv(\phi \bfV) =
    \phi - (\phi/2) \rmdiv(\bfV) - (1/2) \bfV.\nabla \phi
    \geq \phi - (\phi/2) \rmdiv(\bfV) > c \phi.
    $$
    Since
    $\epsilon u_\epsilon \nabla u_\epsilon.\nabla \phi \rightarrow 0$
    in $\Lone$, it follows that $\{u_\epsilon\}_{\epsilon > 0}$
    converges strongly to $u$ and $\sqrt{\epsilon} \nabla u_\epsilon$
    to zero in the $\phi$-weighted $\Ltwo$ space.
    Convergence of $\{u_\epsilon\}_{\epsilon > 0}$ in
    $\Ltwo$ then follows from Lemma~\ref{lem:strong-L2-conv}. 
\end{proof}

\subsection{Piecewise Constant HDG Scheme}\label{sec:pw-const-scheme}
The scheme analyzed in this section uses a hybridized mixed
discretization of the diffusion operator and the classical upwind
discontinuous Galerkin discretization of the transport operator.
The discrete spaces for the piecewise constant HDG scheme are
\begin{equation}\label{eqn:hdg-space}
\calU_h = \calP_0(\calT_h) \times \calP_0(\delTh),
\qquad \text{ and } \qquad
U_h = \{(u_h,\uhat_h) \in \calU_h \sst \uhat_h|_{\partial\Omega} = 0 \}.
\end{equation}
The DG approximation of the transport equation uses
the upwind value on the element interfaces.
\begin{definition}\label{def:upwind}
  Let $k \in \delTh$ be a facet and $u_h \in \calP_0(\calT_h)$.
  \begin{itemize}
  \item If $k = K_1 \cap K_2 \in \delTho$ with $K_1, K_2 \in \calT_h$
    is an interior facet, and $\bfn$ is the outward pointing normal of
    $K_1$, then
    $$
    u^{up}_h|_k = 
    \begin{cases}
      u_h|_{K_1},& \quad \bfV.\bfn \geq 0 \\
      u_h|_{K_2},& \quad \bfV.\bfn < 0
    \end{cases}.
    $$

  \item If $k \subset \partial K \subset \Gammaout$ for some
    $K \in \calT_h$, then $u^{up}_h|_k = u_h|_K$.

  \item If $k \subset \Gammain$ then $u^{up}_h|_k = 0$
    (or more generally $u_0$ for nonhomogeneous inflow boundary data).
  \end{itemize}
  Upwind values on the characteristic boundary do not enter into
  the algorithm, so for definiteness set them to be zero.
\end{definition}

Note that in the mixed formulation, specification of inflow boundary
data for the transport terms enters as $u^{up}|_\Gammain = 0$, while
the Dirichlet data for the diffusion terms appears as
$\uhat_h|_{\partial \Omega} = 0$.

The stabilization term appearing in the HDG method involves the inner
product
$$
\langle u, v\rangle = \sum_{K\in \calT_h}\int_{\partial K} u v,
\qquad
u, v \in \Hone \oplus \calP_0(\delTh).
$$

The degrees of freedom of $\calP_0(\delTh)$ coincide with the degrees
of freedom of the first-order Crouzeix--Raviart space
$CR(\calT_h)$, and gradients appearing in the piecewise constant HDG
scheme coincide with those in the CR space. The latter are characterized as
follows.

\begin{definition}\label{def:cr-lifting}
  The Crouzeix--Raviart lifting operator
  $\mathcal{L}_h: \calP_0(\delTh)\to \calP_1(\calT_h)$ is defined
  by
  $$
  \langle \mathcal{L}_h(\hat{v}_h), \hat{w}_h\rangle = \langle
  \hat{v}_h, \hat{w}_h\rangle \text{ for all } \hat{w}_h\in
  \calP_0(\delTh).
  $$
\end{definition}

Using this lifting and upwind value, the piecewise constant scheme
for~\eqref{eqn:strong0} becomes
\begin{equation}  \label{eqn:hdg0}
(u_h,\uhat_h) \in U_h, \qquad 
a_{h\epsilon}
\left((u_{h\epsilon},\uhat_{h\epsilon}),(v_h,\vhat_h) \vph\right)
= (f,v_h),
\qquad (v_h,\vhat_h) \in U_h,
\end{equation}
where the bilinear form $a_{h\epsilon}: U_h \times U_h \rightarrow \Re$ is
\begin{eqnarray}\label{eqn:hdg-bilin-form}
a_{h\epsilon}
\left((u_{h\epsilon},\uhat_{h\epsilon}),(v_h,\vhat_h) \vph\right) &=&
  (u_{h\epsilon},v_h)
  + \epsilon \langle u_{h\epsilon} - \uhat_{h\epsilon}, v_h - \vhat_h \rangle
  + \epsilon \left(\nabla\calL_h(\uhat_{h\epsilon}),
    \nabla \calL_h(\vhat_h) \vph\right) \nonumber\\
&&  - \left(u_{h\epsilon},  v_h \rmdiv(\bfV) \vph\right)
  + \langle u^{up}_{h\epsilon},  v_h (\bfV.\bfn) \rangle.
\end{eqnarray}

The following assumptions will enable construction of discrete
analogs of the barrier functions used in the continuous setting.

\begin{assumption}\label{as:dirichlet-bc-discrete}
  Assume $\Omega \subset \Re^d$ is a Lipschitz domain and,
  \begin{itemize}
  \item $\{\mathcal T_h\}_{h>0}$ is a regular family of
    triangulations for which every simplex is acute; that is, every
    internal dihedral angle of every simplex $K\in\mathcal T_h$ is at
    most $\pi/2$.

  \item There exists a barrier $\psi$ on a subdomain $\omega$ in the
    sense of Definition~\ref{def:barrier} with the bound on the
    forcing term strengthened to $\psi + \bfV.\nabla \psi \geq |f|$ in
    the sense of distributions on $\omega$.

  \item $\omega$ is fitted by $\omega \cap \calT_h$.
\end{itemize}
\end{assumption}

Under these assumptions, piecewise constant HDG approximations
of the convection-diffusion equation converge strongly in $\Ltwo$.

\begin{theorem}\label{thm:dirichlet-conv-hdg}
  Let Assumption~\ref{as:dirichlet-bc} and Assumption~\ref{as:dirichlet-bc-discrete} hold, and set
  $\phi(x) = \min(\mathrm{dist}(x,\Gammaout),\delta)$. Let
  $u_{h\epsilon}$ be the HDG solution to~\eqref{eqn:hdg0}. As
  $(h, \epsilon)\to(0,0)$, $u_{h\epsilon}\to u$ strongly in $\Ltwo$, where
  $u$ is the (renormalized) solution of the transport equation
  $$
  u + \bfV.\nabla u = f, \quad \text{ in } \calD'(\Omega),
  \quad \text{ with trace } \quad \gamma(u)|_\Gammain = 0.
  $$
  In addition, $\sqrt{\epsilon} \ltwo{\phi^{1/2}
    \nabla \calL_h(\uhat_{h\epsilon})} \rightarrow 0$
  and $\sum_{K \in \calT_h} \int_{\partial K} \epsilon \phi (u_{h\epsilon}
  - \uhat_{h\epsilon})^2\rightarrow 0$ as $(h,\epsilon)\to(0,0)$.
\end{theorem}

\subsubsection{Proof of Theorem \texorpdfstring{\ref{thm:dirichlet-conv-hdg}}{\getrefnumber{thm:dirichlet-conv-hdg}}}
The proof involves many routine, but long, estimates for
discontinuous functions on triangulations. To keep the arguments
transparent, the essential properties required for the proof of the
main theorem are stated next, and their
proofs are provided in the Appendix.

\begin{lemma} \label{lem:dg0}
  Let $\{\calT_h\}_{h > 0}$ be a regular family of meshes for which
  each simplex is acute,
  $$
  \calU_h = \calP_0(\calT_h) \times \calP_0(\delTh)
  \quad \text{ and } \quad
  U_h = \{(u_h,\uhat_h) \sst \uhat_h|_{\partial \Omega} = 0\},
  $$
  and let $a_{h\epsilon}: \calU_h \times \calU_h \rightarrow \Re$ be the
  bilinear form defined in~\eqref{eqn:hdg-bilin-form}. Define
  $\iota_h: \Hone \rightarrow \calU_h$ by
  $\iota_h(v) = (P_hv, \hat{P}_hv)$.
  
  The following hold.
  \begin{enumerate}
  \item {\em Coercivity:} If $(u_h,\uhat_h) \in \calU_h$, then
    \begin{eqnarray*}
      a_{h\epsilon}
      \left((u_h,\uhat_h), (u_h, \uhat_h) \vph\right)
      &=& \int_\Omega (1 - (1/2) \rmdiv(\bfV)) u_h^2
          + \epsilon \langle u_h - \uhat_h, u_h - \uhat_h \rangle^2 
          + \epsilon \ltwo{\nabla \calL_h(\uhat_h)}^2 \\ 
      && + \half \ltwoin{u_h}^2
         + \half \ltwoout{u_h}^2 
         + \half \ltwoDelTho{[u_h]}^2 \\
      &\equiv& \norm{(u_h,\uhat_h)}^2_{h\epsilon}.
    \end{eqnarray*}
    
  \item {\em Monotonicity:} The scheme~\eqref{eqn:hdg0} is monotone in
    the sense that if $(u_h,\uhat_h) \in \calU_h$, then
    \begin{equation}\label{eqn:hdg-monotonicity}
    a_{h\epsilon}
    \left((u_h,\uhat_h), (u_h^+, \uhat_h^+) \vph\right)
    \geq  \norm{(u_h^+,\uhat_h^+)}_{h\epsilon}^2.
    \end{equation}

  \item {\em Consistency:} If $v \in \calD(\Omega)$ and
    $\iota_h(v) = (v_h, \vhat_h)$ is the interpolant of $v$, then
    \begin{multline} \label{eqn:hdg-consistency}
      \ltwo{v - v_h} 
      + \ltwo{\nabla (v - \calL_h(\vhat_h))}
      + \langle v_h - \vhat_h, v_h - \vhat_h \rangle \\
      + \left(\sum_K \int_{\partial K} (v_h - v)^2 |\bfV.\bfn| \right)^{1/2}
      \leq C \left(\hone{\bfV} \vph\right)
      \htwo{v} \sqrt{h}. 
    \end{multline}
  
  \item {\em Inflow Barrier:} If $\psi \in \Wonefour$ vanishes on
    $\Gammain$ and $(\psi_h, \psihat_h) = \iota_h(\psi)$, then
    \begin{equation}\label{eqn:hdg-inflow}
    a_{h\epsilon}\big((\psi_h,\psihat_h), (v_h, \vhat_h) \big)
    = \int_\Omega
    \left(\psi + \bfV.\nabla \psi \vph\right) v_h
    + \epsilon \nabla \calL_h(\psihat_h).\nabla \calL_h(\vhat_h)
    + R^\psi_h,
    \end{equation}
    where the remainder $R^\psi_h$ is bounded as
    \begin{equation}\label{eqn:hdg-inflow-remainder}
    |R^\psi_h| \leq C \left(\hone{\bfV}, \linf{\rmdiv(\bfV)} \vph\right)
    \lfour{\nabla \psi}
    \norm{(v_h,\vhat_h)}_{h\epsilon} \sqrt{h}.
    \end{equation}
    
  \item {\em Outflow Cutoff:} If
    $\phi(x) = \min(\mathrm{dist}(x, \Gammaout), \delta)$ and
    $\iota_h(\phi) = (\phi_h, \phihat_h) \in U_h$, then
    \begin{eqnarray}
      \lefteqn{
      a_{h\epsilon}\big((u_h,\uhat_h), (\phi_h u_h, \phihat_h \uhat_h) \big)
      } \nonumber\\
    &=& \int_\Omega \left(\phi - (1/2) \rmdiv(\phi \bfV) \vph\right) u_h^2
    + \epsilon \phi |\nabla \calL_h(\uhat_h)|^2
    + \epsilon u_h \nabla \calL_h(\uhat_h) . \nabla \calL_h(\phi)  \label{eqn:hdg-outflow-cutoff}\\
      &&  + \sum_{k \in \delTho} \int_k\half \phi [u_h]^2 |\bfV.\bfn|
         + \int_\Gammain \!\! \half\phi \, u_h^2 |\bfV.\bfn|
         + \epsilon \sum_K \int_{\partial K} \phi (u_h - \uhat_h)^2 
     + R^\phi_h,\nonumber
  \end{eqnarray}
  where the remainder $R^\phi_h$ is bounded as
  $$
  |R^\phi_h| \leq C\hone{\bfV}^{1/2} \vph
  \linf{\nabla \phi}
  \norm{(u_h, \uhat_h)}_{h\epsilon}^2 h^{1-d/4}(\log(h^{-1}))^{1/4}.
  $$
\end{enumerate}
\end{lemma}

The following elementary lemma is proved in the Appendix.

\begin{lemma}\label{lem:convex-cone}
  Let $\Omega \subset \Re^d$, and let $\{u_h\}_{h>0} \subset \Ltwo$
  converge weakly to a limit $u$. If $\psi \in \Ltwo$ and
  $(u_h \mp \psi)^\pm$ both converge strongly to zero, then
  $-\psi \leq u \leq \psi$.
\end{lemma}

Granted these lemmas, the proof of Theorem~\ref{thm:dirichlet-conv-hdg}
now parallels the proof of Theorem~\ref{thm:dirichlet-conv-cont}.

{\em1. Weak Limit:} Since $\linf{\rmdiv(\bfV)} < 2$, there is a
constant $C>0$ such that
$\ltwo{u_{h\epsilon}} \leq C
\norm{(u_{h\epsilon},\uhat_{h\epsilon})}_{h\epsilon}$. In this case,
  $$
  \norm{(u_{h\epsilon},\uhat_{h\epsilon})}_{h\epsilon}^2
  = a_{h\epsilon}\big((u_{h\epsilon},\uhat_{h\epsilon}),
  (u_{h\epsilon}, \uhat_{h\epsilon}) \big)
  = (f,u_{h\epsilon}) \leq C \ltwo{f}
  \norm{(u_{h\epsilon},\uhat_{h\epsilon})}_{h\epsilon}.
  $$
  Hence
  $\{\norm{(u_{h\epsilon},\uhat_{h\epsilon})}_{h\epsilon}\}_{h,\epsilon
    > 0}$ is bounded, and upon passing to a subsequence we obtain a
  weakly converging sequence $u_{h\epsilon} \weak u$ in $\Ltwo$.

  When $f \in \Linf$, additional bounds are available for solutions of~\eqref{eqn:hdg0}. If $k \in \Re$, then
  $$
  a_{h\epsilon}\big((u_{h\epsilon}-k,\uhat_{h\epsilon}-k),
  (v_h, \vhat_h) \big)
  = (f-k, v_h),
  \qquad (v_h, \vhat_h) \in U_h.
  $$
  If $k \geq 0$, then $(v_h, \vhat_h) \equiv \left( (u_h-k)^+, (\uhat-k)^+
    \vph\right) \in U_h$ since $\uhat_{h\epsilon} - k \leq
  0$ on the boundary. If $k \geq \linf{f}$, then
  $$
  \norm{((u_{h\epsilon}-k)^+, (\uhat_{h\epsilon}-k)^+)}_{h\epsilon}^2
  \leq a_{h\epsilon}\big((u_{h\epsilon}-k,\uhat_{h\epsilon}-k),
  ((u_{h\epsilon}-k)^+, (\uhat_{h\epsilon}-k)^+) \big)
  \leq (f-k, (u_{h\epsilon} - k)^+) \leq 0,
  $$
  so $u_{h\epsilon} \leq \linf{f}$. The same argument shows
  $u_{h\epsilon} \geq -\linf{f}$, whence
  $\linf{u_{h\epsilon}} \leq \linf{f}$. In particular,
  $u_{h\epsilon} \weak u$ in $\Lp$ for $1 \leq p < \infty$.
  
  {\em2. Limit Equation:} Fix $v \in \calD(\Omega)$ and let
  $(v_h, \vhat_h) = \iota_h(v) \in U_h$ be the projection in Lemma~\ref{lem:dg0}. Write the weak statement~\eqref{eqn:hdg0} satisfied
  by $(u_{h\epsilon}, \uhat_{h\epsilon})$ as
  $$
  (u_{h\epsilon},v_h)
  + \sum_{K \in \calT_h} \int_{\partial K} (u^{up}_{h\epsilon}-u_{h\epsilon})
  (\bfV.\bfn) v_h
  = (f,v_h)
  - \epsilon \langle u_{h\epsilon} - \uhat_{h\epsilon}, v_h - \vhat_h \rangle
  - \epsilon \left(\nabla \calL_h(\uhat_{h\epsilon}),
    \nabla \calL_h(\vhat_h) \vph\right),
  $$
  and rewrite the second term on the left as
  \begin{eqnarray*}
    \sum_{K \in \calT_h} \int_{\partial K} (u^{up}_{h\epsilon}-u_{h\epsilon})
    (\bfV.\bfn) v_h
    &=& \sum_{K \in \calT_h} \int_{\partial K}
        (u^{up}_{h\epsilon}-u_{h\epsilon})
        (\bfV.\bfn) (v_h - v)
        - u_{h\epsilon} (\bfV.\bfn) v \\
    &=& \sum_{K \in \calT_h} \int_{\partial K}
        (u^{up}_{h\epsilon}-u_{h\epsilon})
        (\bfV.\bfn) (v_h - v)
        - \int_{K} u_{h\epsilon} \left(\rmdiv (v \bfV) \vph\right).
  \end{eqnarray*}
  The first line follows since the quantity
  $u^{up}_{h\epsilon} v \bfV$ is single-valued on each interface, in
  which case the dot products with each of the normals are ``equal and
  opposite''. The second line follows from the divergence theorem and
  the fact that $u_{h\epsilon}$ is constant on $K$. Therefore,
  \begin{multline*}
    \left(u_{h\epsilon} ,v - \rmdiv(v \bfV) \vph\right)
    = (f,v_h) \\
    - \epsilon \langle u_{h\epsilon} - \uhat_{h\epsilon}, v_h - \vhat_h \rangle
    - \epsilon \left(\nabla \calL_h(\uhat_{h\epsilon}),
      \nabla \calL_h(\vhat_h) \vph\right)
    - \sum_{K \in \calT_h} \int_{\partial K} (u^{up}_{h\epsilon}-u_{h\epsilon})
    (\bfV.\bfn) (v_h - v).
  \end{multline*}
  Since
  $\{\norm{(u_{h\epsilon},\uhat_{h\epsilon})}_{h\epsilon}\}_{h,\epsilon}$
  is bounded, the consistency bound~\eqref{eqn:hdg-consistency} in
  Lemma~\ref{lem:dg0} shows that each of the terms in the second line
  vanishes as $(h,\epsilon) \rightarrow (0,0)$. Hence the weak limit
  satisfies $u + \bfV.\nabla u = f$ in $\calD'(\Omega)$.

  {\em3. Boundary Values:} Let
  $\iota_h(\psi) = (\psi_h, \psihat_h) \in \calU_h$ be the projection
  of the barrier function. If $(v_h,\vhat_h) \in U_h$ and
  $v_h \geq 0$, it follows from~\eqref{eqn:hdg-inflow} in Lemma~\ref{lem:dg0} that
  \begin{eqnarray*}
    a_{h\epsilon}\big((u_{h\epsilon}-\psi_h, \uhat_{h\epsilon}
    - \psihat_h), (v_h, \vhat_h) \big)
    &=& (f,v_h) 
    - a_{h\epsilon}\big((\psi_h,\psihat_h), (v_h, \vhat_h) \big) \\
    &=& \left(f - \psi - \bfV.\nabla \psi, v_h \vph\right)
        - \epsilon (\nabla \calL_h(\psihat_h), \nabla \calL_h(\vhat_h))
        + R^\psi_h \\
    &\leq& - \epsilon (\nabla \calL_h(\psihat_h), \nabla \calL_h(\vhat_h))
           + R^\psi_h \\
    &\leq& C \sqrt{\epsilon} \ltwo{\nabla \psi}
           \norm{(v_h^+, \vhat_h^+)}_{h\epsilon}
           + R^\psi_h,
  \end{eqnarray*}
  where the first inequality follows since
  $\psi + \bfV.\nabla \psi \geq |f|$ and $\nabla \psi$ has support in
  $\omega$. Set
  $(v_h, \vhat_h) = ((u_{h\epsilon} - \psi_{h\epsilon})^+,
  (\uhat_{h\epsilon} - \psihat_{h\epsilon})^+)$. Monotonicity~\eqref{eqn:hdg-monotonicity} and the bound~\eqref{eqn:hdg-inflow-remainder} on $R^\psi_h$ from Lemma~\ref{lem:dg0} then show
  \begin{eqnarray*}
    \norm{((u_{h\epsilon}-\psi_h)^+,
    (\uhat_{h\epsilon} - \psihat_h)^+)}_{h\epsilon}^2
    &\leq& a_{h\epsilon}\big(
           (u_{h\epsilon}-\psi_h, \uhat_{h\epsilon} - \psihat_h),
           (v_h, \vhat_h) \big) \\
    &\leq& C \norm{((u_{h\epsilon}-\psi_h)^+,
           (\uhat_{h\epsilon} - \psihat_h)^+)}_{h\epsilon}
           ( \sqrt{\epsilon} + \sqrt{h} ).
  \end{eqnarray*}
  Hence
  $$
  \norm{(u_{h\epsilon} - \psi_h)^+}_{\Ltwo}^2
  \leq C \norm{((u_{h\epsilon} - \psi_h)^+,
    (\uhat_{h\epsilon} - \psihat_h)^+)}^2_{h\epsilon}
  \leq C ( \sqrt{\epsilon} + \sqrt{h} )
  $$
  and similarly
  $$
  \norm{(u_{h\epsilon} +\psi_h)^-}_{\Ltwo}^2
  \leq C\norm{((u_{h\epsilon} + \psi_h)^-,
    (\uhat_{h\epsilon} + \psihat_h)^-)}^2_{h\epsilon}
  \leq C ( \sqrt{\epsilon} + \sqrt{h} ).
  $$
  Lemma~\ref{lem:convex-cone} then shows $-\psi \leq u \leq \psi$, and
  since $\psi$ vanishes on $\Gammain$, it
  follows that the trace $\gamma(u)$ also vanishes on the inflow
  boundary. In particular, $u$ is the solution of
  $u + \bfV.\nabla u = f$ in $\calD'(\Omega)$ with
  $\gamma(u)|_{\Gammain} = 0$, and since the solution is unique, it
  follows that the whole sequence $\{u_{h\epsilon}\}_{h,\epsilon > 0}$
  converges weakly to $u$.

  {\em4. Strong Convergence:} To establish strong convergence, let
  $\phi(x) = \min(\mathrm{dist}(x, \Gammaout), \delta)$.
  Taking the limit superior of~\eqref{eqn:hdg-outflow-cutoff} as
  $h,\epsilon\to0$ yields
  \begin{align*}
    \limsup_{(h,\epsilon) \rightarrow (0,0)}
      \Big\{\int_\Omega&
      \left(\phi - (1/2) \rmdiv(\phi \bfV) \vph\right) u_{h\epsilon}^2
         + \phi \epsilon |\nabla \calL_h(\uhat_{h\epsilon})|^2
         +u_{h\epsilon} \epsilon \nabla \calL_h(\uhat_{h\epsilon}) .
          \nabla \calL_h(\phi)\\
      &+ \sum_{K \in \calT_h} \int_{\partial K}
         \epsilon \phihat_h (u_{h\epsilon} - \uhat_{h\epsilon})^2
         + \sum_{k \in \delTho} \int_k
         (\phi/2) [u_{h\epsilon}]^2 |\bfV.\bfn|
        + \int_\Gammain \!\! \half \phi \, u_{h\epsilon}^2 |\bfV.\bfn|
        + R^\phi_h \Big\} \\
      &= \limsup_{(h,\epsilon) \rightarrow (0,0)}
          a_{h\epsilon}\big((u_{h\epsilon},\uhat_{h\epsilon}),
          (\phi_h u_{h\epsilon}, \phihat_h \uhat_{h\epsilon}) \big) \\
          &= \limsup_{(h,\epsilon) \rightarrow (0,0)}
          (f, \phi_h u_{h\epsilon}) = (f, \phi u) 
      = \int_\Omega \left(\phi - (1/2) \rmdiv(\phi \bfV) \vph\right) u^2,
    \end{align*}
    where the last equality follows from~\eqref{eqn:test-by-u-phi}.
    The residual estimate~\eqref{eqn:hdg-inflow-remainder} on
    $R^\psi_h$ implies that $R^\psi_h\to0$ as $h,\epsilon \to0$, and
    the elementary estimate
    $$
    \epsilon \left(u_{h\epsilon} \nabla \calL_h(\uhat_{h\epsilon}),
      \nabla \calL_h(\phihat_h) \vph\right)
    \leq C \norm{(u_{h\epsilon}, \uhat_{h\epsilon})}_{h\epsilon}
    \ltwo{\nabla \phi} \sqrt{\epsilon}
    $$
    also implies that
    $\epsilon (u_{h\epsilon} \nabla \calL_h(\uhat_{h\epsilon}), \nabla
    \calL_h(\phihat_h) \vph)$ vanishes in the limit. Strong
    convergence of $u_{h\epsilon}$, the jumps, and the trace on
    $\Gammain$, in $\phi$-weighted $L^2$ norms, now
    follows. Additionally,
    $\int_\Omega \epsilon \phi |\nabla \calL_h(\uhat_{h\epsilon})|^2$
    and
    $\sum_{K \in \calT_h} \int_{\partial K} \epsilon \phi
    (u_{h\epsilon} - \uhat_{h\epsilon})^2$ converge to zero.
    Convergence of $u_{h\epsilon}\to u$ strongly in (unweighted)
    $\Ltwo$ as $(h,\epsilon)\to (0,0)$ follows from Lemma~\ref{lem:strong-L2-conv}.

\subsection{Computational Example}

We present a computational example using NGSolve~\cite{schoberl2014c++} to illustrate the results of our
main theorem. We let $\Omega = (0,1)^2$ and consider
the velocity field 
\begin{equation*}
  \bfV(\bfx) = (1-\rho(\bfx))\bfe_2 
  + \rho(\bfx) r^{1/10} \big( \sin(\theta), -\cos(\theta)\big)^\top,
\end{equation*}
where $r, \theta$ are polar coordinates using $(0.5,0.5)$ as the
center, and $\rho$ is a $C^1(\Omegabar)$ cutoff function with support
in $B_{1/4}(0.5,0.5)$. Then $\bfV \in \Hone^d\cap \Linf^d$,
$\mathrm{div}(\bfV) \in \Linf$, and
$\bfV \not\in W^{1,\infty}(\Omega)$. The solution is manufactured to be
$u(\bfx) = x_2-4x_2 \left(x_1-1/2\right)^2$ by setting
$f = u +\bfV.\nabla u$.

\begin{figure}
\begin{center}
\includegraphics[width=.3\textwidth]{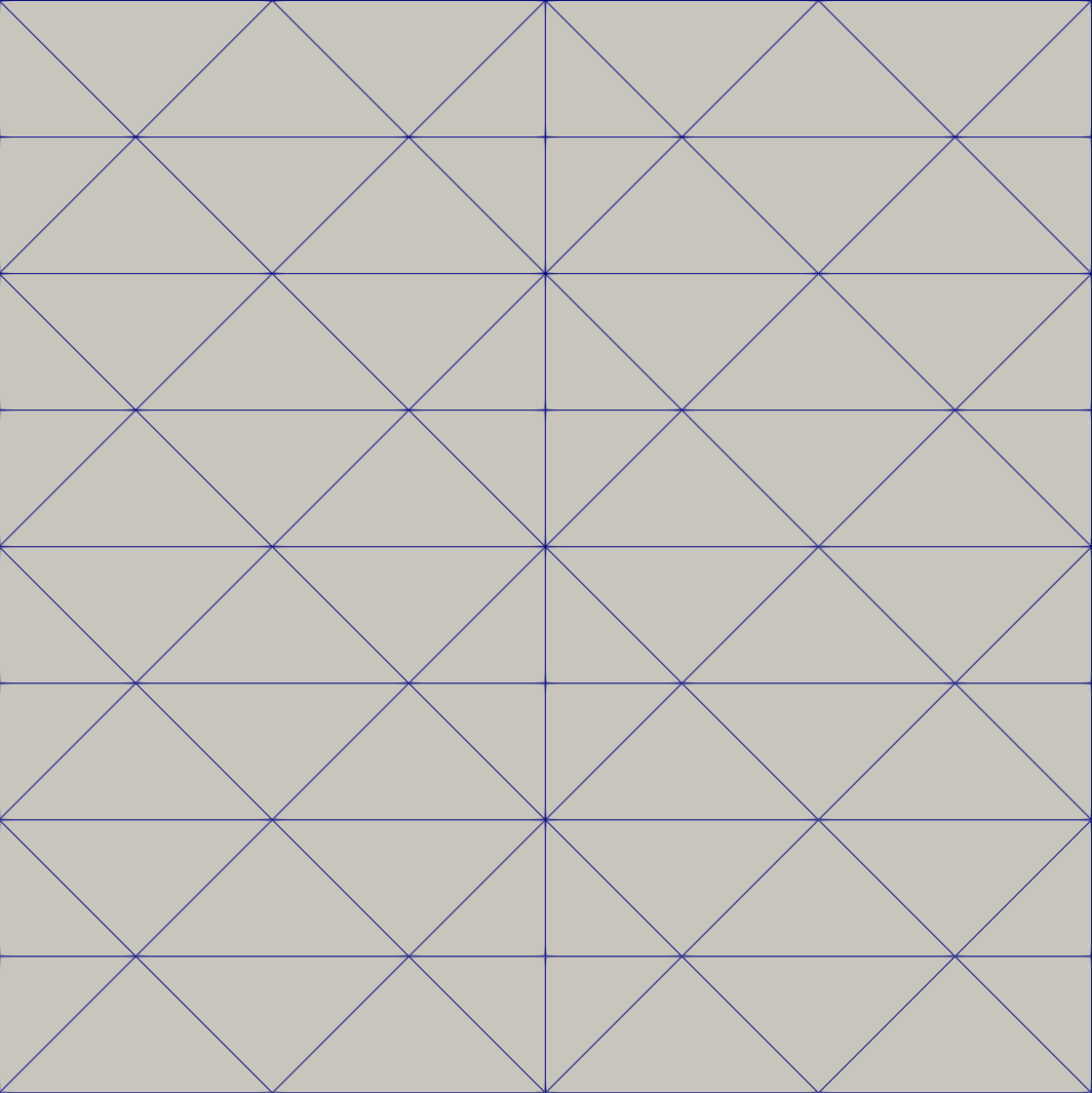}
\includegraphics[width=.3\textwidth]{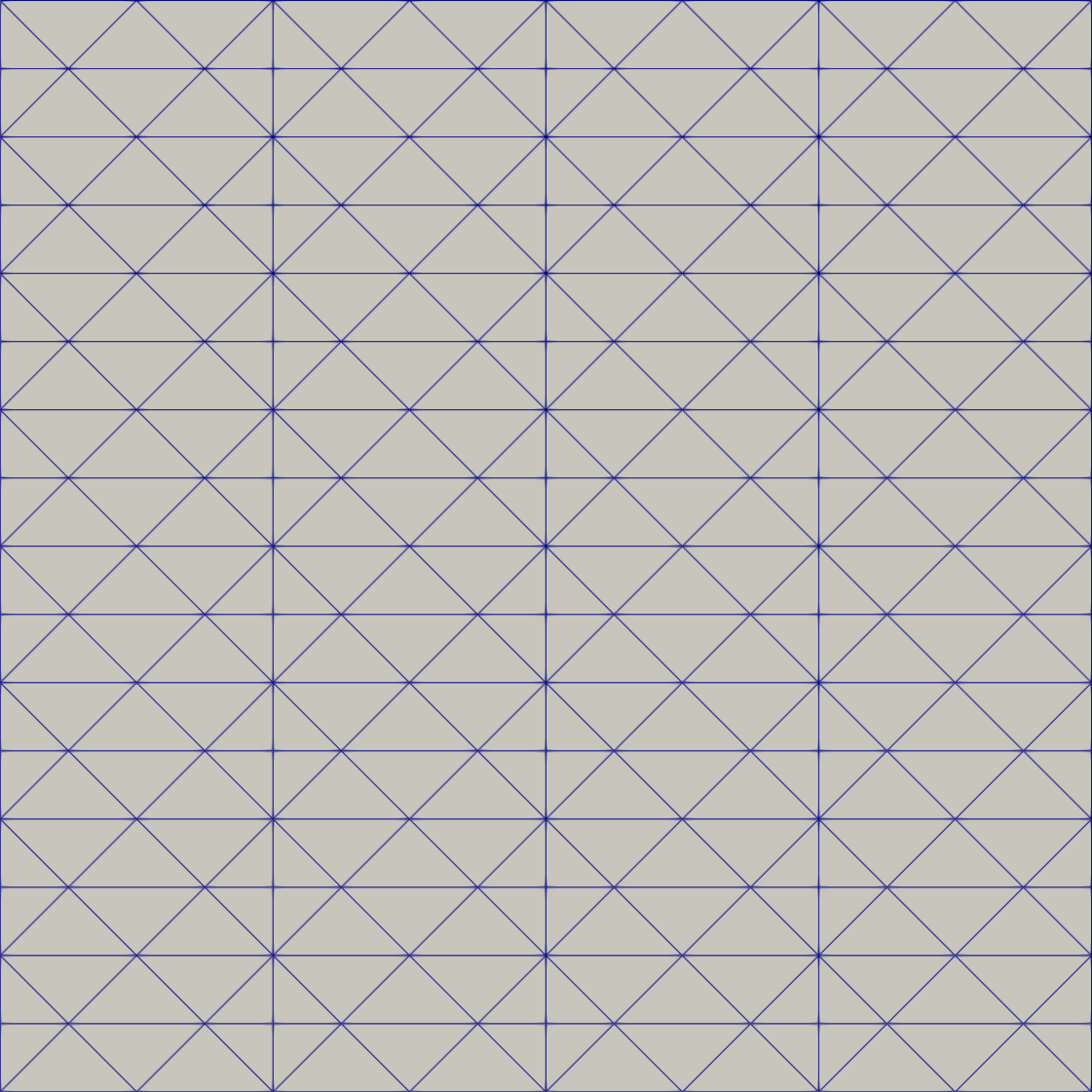}
\includegraphics[width=.3\textwidth]{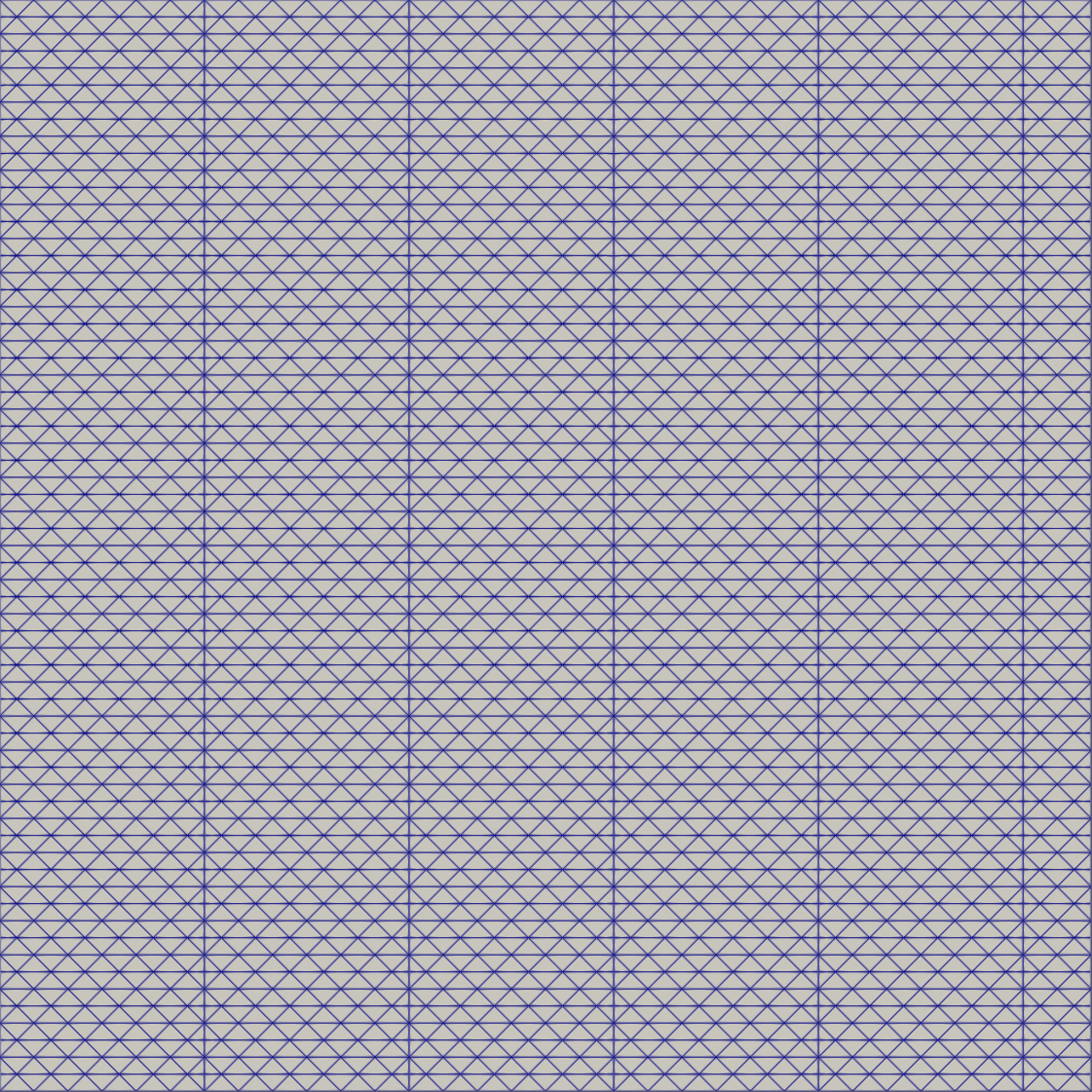}
\end{center}
\caption{Meshes constructed for the example
  following~\cite{peterson1991note}. Left-to-right, $h=1/4,1/8,1/32$, with $m=1,3,5$ vertical lines.}
\label{fig:peterson-mesh}
\end{figure}

To observe the usual $O(\sqrt h)$ rate~\cite{johnson1986analysis} for
$\Vert u - u_{h\epsilon}\Vert_\Ltwo$, we use the sequence of meshes
introduced by Peterson~\cite{peterson1991note}. The domain is meshed
with triangles and additional vertical lines, the number of which
increases as $h \rightarrow 0$.  The solution vanishes on the inflow
boundary, $\Gammain = \{(x,0) \sst 0 \leq x \leq 1\}$, and the
triangles are not aligned with the flow field of $\bfV$ in the
sense of~\cite{cockburn2008optimal}.  Figure~\ref{fig:peterson-mesh}
displays three meshes with different $h$ values and number of vertical
lines.

\begin{figure}
\begin{center}
\includegraphics[height=1.6in]{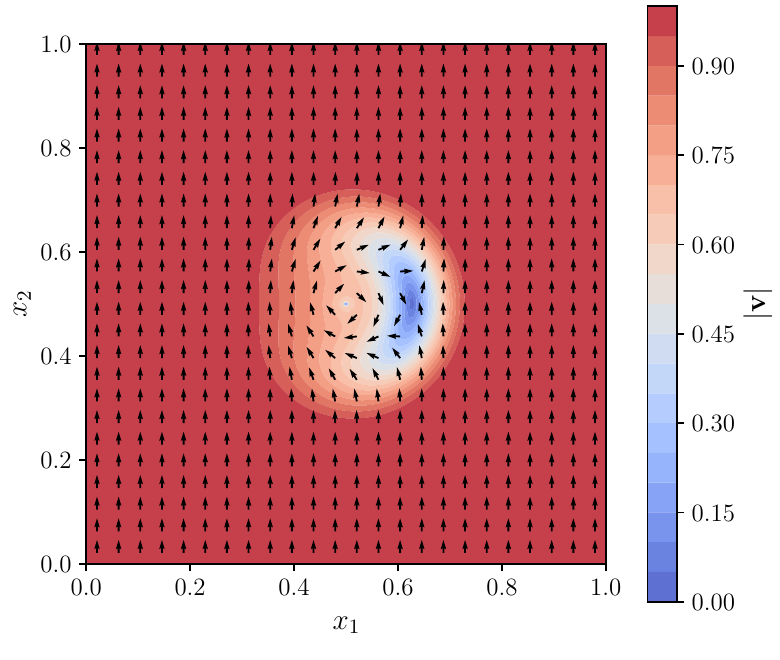}
\includegraphics[height=1.6in]{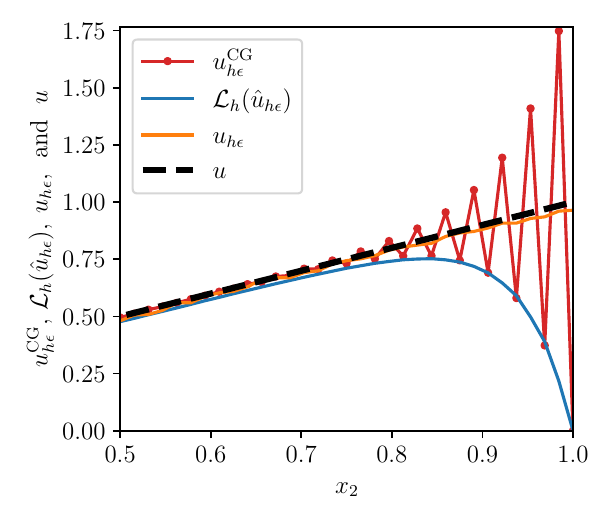}
\includegraphics[height=1.6in]{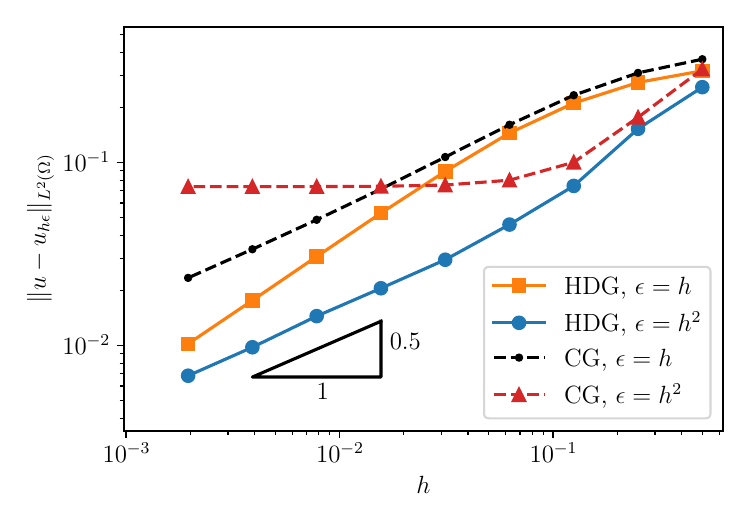}
\end{center}
\caption{(Left) The velocity field $\bfV$. (Middle) Plot of the discrete solutions on the slice $\{x_1 = 1/2\}$ near the outflow boundary for $h = 2^{-5}$ and $\epsilon = h^2$.
(Right) A log-log plot of the $\Ltwo$ errors vs.\ $h$.
}
\label{fig:dirichlet}
\end{figure}

For this experiment, we take $h = 2^{-\ell}$ for $\ell=1,\ldots,9$ and
$m = \texttt{round}(h^{-0.43})$ for the piecewise constant HDG scheme
and a CG scheme, with $\epsilon = h$ and $\epsilon = h^2$ for both
schemes. Figure~\ref{fig:dirichlet} shows the velocity field $\bfV$, a
plot of the discrete solutions on the slice $\{x_1 = 1/2\}$ near the
outflow boundary, and a plot of the $\Ltwo$ errors vs.\ $h$. In
Figure~\ref{fig:dirichlet} (middle), the CG solution,
$u_{h\epsilon}^{\mathrm{CG}}$, oscillates, but the HDG solution $u_{h\epsilon}$
stays close to $u$. The CR reconstruction
$\mathcal{L}_h(\hat{u}_{h\epsilon})$ is overdiffused and is not a good
approximation of $u$. Our theorems do not say anything about
$\hat{u}_{h\epsilon}$. At $x_2 = 0.5$, $\bfV$ is not locally
Lipschitz. Despite this, $u_{h\epsilon}^{\mathrm{CG}}$ is well-behaved
near $x_2 = 0.5$. It is primarily the boundary conditions that cause
the oscillations in
$u_{h\epsilon}^{\mathrm{CG}}$. Figure~\ref{fig:dirichlet} (right)
shows that when $\epsilon = h^2$, the CG solution error stalls,
matching Example~\ref{ex:cg-counter-example}.  The HDG solution for
both scalings of $\epsilon$ appears to converge at the rate
$\sqrt{h}$. This scaling matches the standard error estimate for
classical DG applied to the transport equation.

\section*{Declarations} The authors developed the mathematical ideas (applying the renormalization theory and trace operator of~\cite{boyer2005trace} and the continuous and discrete barrier arguments) and wrote the first draft of the manuscript without the use of AI. 
In the process of revising the manuscript, the first author used GPT-5.6 Sol (OpenAI) and Claude Opus 5 (Anthropic) to assist with literature search, 
reviewing the mathematical arguments, and writing plotting scripts. 
During this process, GPT suggested an idea that led to the proof of Lemma~\ref{lem:strong-L2-conv}, which the authors wrote and verified. Consequently, the authors improved Theorems~\ref{thm:dirichlet-conv-cont} and~\ref{thm:dirichlet-conv-hdg} from a weighted $L^2$ convergence result (weighted with $\phi$) to an unweighted $L^2$ convergence result. 
Claude and GPT pointed out the relevant references \cite{bacuta2026explaining} and \cite{schlichting2018analysis}.
The authors assume responsibility for all content.

%%%%%%%%%
\appendix

\section{Proofs of Lemmas}

This section starts with proofs of the elementary Lemmas~\ref{lem:strong-L2-conv} and~\ref{lem:convex-cone}
used in the convergence proofs.

\begin{proof}[Proof of Lemma~\ref{lem:strong-L2-conv}] Let $\eta > 0$
  and split the integral as
\begin{align*}
  \Vert u - u_h \Vert_{\Ltwo}^2
  &= \int_{\{ \phi \leq \eta\}} (u - u_h)^2
    + \int_{\{ \phi > \eta\}} (u - u_h)^2
    \leq C^2|\{ \phi \leq \eta \}|
    + \frac{1}{\eta} \int_{\{ \phi > \eta\}} \phi (u - u_h)^2,
\end{align*}
where $C = 2\sup_{h>0}\Vert u_h\Vert_{L^\infty(\Omega)} <\infty$.
Passing to the limit as $h\to0$ implies $\limsup_{h\to0} \Vert u - u_h \Vert_{\Ltwo}^2 \leq C^2|\{ \phi \leq \eta \}| $ for all $\eta>0$. By continuity from above of the Lebesgue measure, $\lim_{\eta\to0}|\{ \phi \leq \eta \}| = |\{ \phi = 0 \}| = 0$. Hence, $\lim_{h\to0} \Vert u - u_h \Vert_{\Ltwo} = 0$.
\end{proof}

\begin{proof}[Proof of Lemma~\ref{lem:convex-cone}] It is sufficient to
  show $(u \mp \psi)^\pm = 0$. We first show $(u - \psi)^+=0$. Let
  $w \in \Ltwo$ be nonnegative. Then
  \begin{equation*}
    \int_\Omega (u - \psi) w
    = \lim_{h \rightarrow 0} \int_\Omega (u_h - \psi) w
    \leq  \lim_{h \rightarrow 0} \int_\Omega (u_h - \psi)^+ w \leq 0.
  \end{equation*}
  Setting $w = (u-\psi)^+$ shows $(u-\psi)^+ = 0$, and the same
  argument applied to $-(u+\psi)$ shows $(u+\psi)^- = 0$.
\end{proof}

\subsection{Piecewise Constant Approximation}
This section collects properties of piecewise constant polynomials and
the DG and HDG approximations utilized above. The proof of
Lemma~\ref{lem:dg0} will use the following three elementary
lemmas.

We first recall the following result from~\cite[Lemma~2.4]{jiang2025stability}.

\begin{lemma} \label{lem:interp} Let
  $\{\calT_h\}_{h > 0}$ be a regular family of meshes, let
  $1 \leq p < \infty$, and let
  $\iota_h:\Wonep \rightarrow \calP_0(\calT_h) \times \calP_0(\delTh)$
  be the projection operator defined in Lemma~\ref{lem:dg0}.  Then
  $$
  \lp{v - v_h} \leq C h \lp{\nabla v},
  \quad \text{ and } \quad
  \lp{\nabla\calL_h(v) - \nabla v} \rightarrow 0.
  $$
  In addition,
  $$
  \left(\sum_{K \in \calT_h}
      \int_{\partial K} |v_h - \vhat_h|^p \right)^{1/p}
    \leq C \lp{\nabla v} h^{1-1/p};
    $$
    in particular, if $v \in \Hone$ then
    $\langle v_h - \vhat_h \rangle \leq C \sqrt{h} \ltwo{\nabla v}$, and 
      \begin{align*}
  \ltwo{v-v_h} + \sqrt{h}  \langle v_h - \vhat_h\rangle &\leq C \ltwo{\nabla v} h, \\
  \ltwo{\nabla v - \nabla \calL_h(\vhat_h)} &\leq C \ltwo{D^2 v} h.
  \end{align*}
\end{lemma}

If $K \in \calT_h$ is a simplex, the inequalities
\begin{align*}
  \norm{u}_{L^2(\partial K)}
  &\leq (C/ \sqrt{h}) \norm{u}_{H^1(K)}, \qquad u \in H^1(K), \\
  \norm{u_h}_{L^2(\partial K)}
  &\leq (C/ \sqrt{h}) \norm{u_h}_{L^2(K)}, \qquad u_h \in \calP_k(K),
\end{align*}
where $C$ depends only upon the aspect ratio of $K$, will be used
repeatedly below. The following sharper estimate is required to
bound the cutoff function.

\begin{lemma} \label{lem:trace}
  Let $K \subset \Re^d$ be a simplex, and $u \in H^1(K)$. There exists
  $C > 0$ depending only upon $d$ and the aspect ratio of $K$ such
  that
  $$
  \int_{\partial K} |u| \leq C 
  \left( h^{d-1-d/p}\norm{u}_{L^{p}(K)}
    + h^{d/2}\norm{\nabla u}_{L^2(K)} \vph\right),
  $$
  where $p$ is less than the Sobolev embedding exponent $2^*$, i.e.\
  $1/2^* = \max\{1/2 - 1/d, 0\}$, and $h = \mathrm{diam}(K)$.
  
  In particular, if $\calT_h$ is a regular triangulation of $\Omega \subset
  \Re^d$ with $d \geq 2$, then
  $$
  \max_{K \in \calT_h}   \int_{\partial K} |u|
  \leq \begin{cases}
  C \hone{u} h\sqrt{\log(h^{-1})}, &d = 2 \\
   C \hone{u} h^{d/2},  &d\geq 3
   \end{cases}.
  $$
\end{lemma}

\begin{proof}
  Mapping $u \in H^1(K)$ to and from a parent element with an affine
  transformation establishes the existence of a constant $C > 0$
  such that
  $$
  \int_{\partial K} |u| \leq \frac{C}{h} \int_K |u| + h |\nabla u|.
  $$
  Bounding the integrals on the right by the higher Lebesgue norms
  gives
    \begin{eqnarray*}
    \int_{\partial K} |u| 
    &\leq& C \big( \norm{u}_{L^{p}(K)} |K|^{1-\frac{1}{p}} / h
         + \norm{\nabla u}_{L^2(K)} |K|^{1/2} \big) \\
    &\leq& C \big( h^{d-1-d/p}\norm{u}_{L^{p}(K)} 
           + h^{d/2}\norm{\nabla u}_{L^2(K)} \big)  ,
  \end{eqnarray*}
  where the last line follows since $|K| \leq C h^d$.  For $d=2$, we
  pick $p_h = \log(h^{-1})$. The Sobolev embedding constant for
  $H^1(\Omega)$ into $L^p(\Omega)$ for $d=2$ scales like
  $C_{p_h} = C\sqrt{p_h} = C\sqrt{\log(h^{-1})}$. Hence,
  \begin{equation*}
    \int_{\partial K} |u| 
    \leq C \left(h^{1-\frac{2}{p_h}}\Vert u\Vert_{L^{p_h}(K)}
      + h\Vert \nabla u\Vert_{L^2(K)} \right)
    \leq Ch^{1-\frac{2}{{p_h}}}\sqrt{\log(h^{-1})} \Vert  u\Vert_{H^1(K)}.
  \end{equation*}
  Since
  $h^{-\frac{2}{{p_h}}} = e^{\frac{2}{\log(h^{-1})}\log(h^{-1})} =
  e^2$, the upper bound simplifies to
  \begin{equation*}
    \int_{\partial K} |u|  \leq C h\sqrt{\log(h^{-1})}\Vert  u\Vert_{H^1(K)}
    = C h^{d/2}\sqrt{\log(h^{-1})}\Vert  u\Vert_{H^1(K)}.
  \end{equation*}
  For $d\geq 3$, we may pick $p = 2^*$, where $2^*$ is the Sobolev
  embedding exponent. Hence $\frac{-d}{p} = 1 - \frac{d}{2}$, and
  \begin{equation*}
    \int_{\partial K} |u| 
    \leq C h^{d/2}\norm{u}_{H^1(K)}.
  \end{equation*}
  The proof of the global inequality follows from
  $\norm{u}_{H^1(K)} \leq \norm{u}_{H^1(\Omega)} $ and from taking the
  maximum over all $K$.
\end{proof}

\begin{lemma}\label{lem:interior-edges-jumps}
  Let $u_h\in \calP_0(\mathcal{T}_h)$, and let $\phi\in \Hone$. The
  following identity and inequality hold:
  \begin{align}
    \sum_{K \in \calT_h} \int_{\partial K \setminus \partial\Omega}
    (u^{up}_h - u_h/2) u_h \phi\, \bfV.\bfn \;
    &= \; \frac{1}{2} \sum_{k \in \delTho}  \int_k [u_h]^2
      \phi |\bfV.\bfn|, \label{eqn:jumps-equality}\\
    \sum_{K \in \calT_h} \int_{\partial K \setminus \partial\Omega}
    (u^{up}_h - u_h/2) u_h^+ \, \bfV.\bfn \;
    &\geq \; \frac{1}{2} \sum_{k \in \delTho}  \int_k [u_h]^2 |\bfV.\bfn|.
      \label{eqn:jumps-inequality}
\end{align}
\end{lemma}
\begin{proof}
  Let $v_h\in \calP_0(\mathcal{T}_h)$ and consider a single facet
  $k\in \delTho$ shared by $K^{up}$ and $K$, with $\mathbf{N}$ the
  uniquely defined normal satisfying $\bfV.\mathbf{N}\geq 0$, as shown
  in Figure~\ref{fig:jumps}. Then $u_h = u^{up}_h$ on $K^{up}$, and
  \begin{align}
    \int_{\partial K^{up} \cap k}
    (u^{up}_h - u_h/2) v_h^{up} \phi \bfV.\bfn
    +  \int_{\partial K \cap k}
    & (u^{up}_h - u_h/2) v_h \phi\bfV.\bfn
      \nonumber\\
    &  = \int_k \bigg( (u^{up}_h/2) v^{up}_h
      - (u^{up}_h - u_h/2) v_h\bigg) \phi|\bfV.\bfN|.
      \label{eqn:jumps-intermediate}
  \end{align}
  In the case of $v_h = u_h$, the identity
  $(a/2)a - (a - b/2)b = a^2/2 -ab +b^2/2 = (1/2) (a-b)^2$ gives
  \begin{align*}
    \int_{\partial K^{up} \cap k} (u^{up}_h - u_h/2) u_h \phi \bfV.\bfn
    +  \int_{\partial K \cap k} (u^{up}_h - u_h/2) u_h \phi\bfV.\bfn
    = \frac{1}{2} \int_k [u_h]^2 \phi|\bfV.\bfn|.
  \end{align*}
  Summing over $k \in \delTho$ proves~\eqref{eqn:jumps-equality}.

  In the case $\phi = 1$ and $v_h = u_h^+$, the integrand on the right-hand side
  of~\eqref{eqn:jumps-intermediate} reads
  $$
  (u^{up}_h/2) v^{up}_h - (u^{up}_h - u_h/2) v_h
  = \frac{1}{2}\bigg(u^{up}_h(u^{up}_h)^+ - 2u^{up}_h u_h^+ + u_hu_h^+\bigg).
  $$
The inequality~\eqref{eqn:jumps-inequality} now follows from
the elementary inequality
  $$
  a a^+ - 2 a b^+ + b b^+ 
  = (a^+)^2 - 2 a b^+ + (b^+)^2 
  = (a^+ - b^+)^2 + 2 (a^+ - a) b^+ \geq (a^+ - b^+)^2,
  $$
and summing over $k \in \delTho$.
\end{proof}

\subsection{Proof of Lemma \ref{lem:dg0}}

\begin{figure}
  \begin{center}
    \scalebox{0.9}{\input{jumps}}
    \caption{Computing the jump terms; upwind and downwind values.}
    \label{fig:jumps}
  \end{center}
\end{figure}
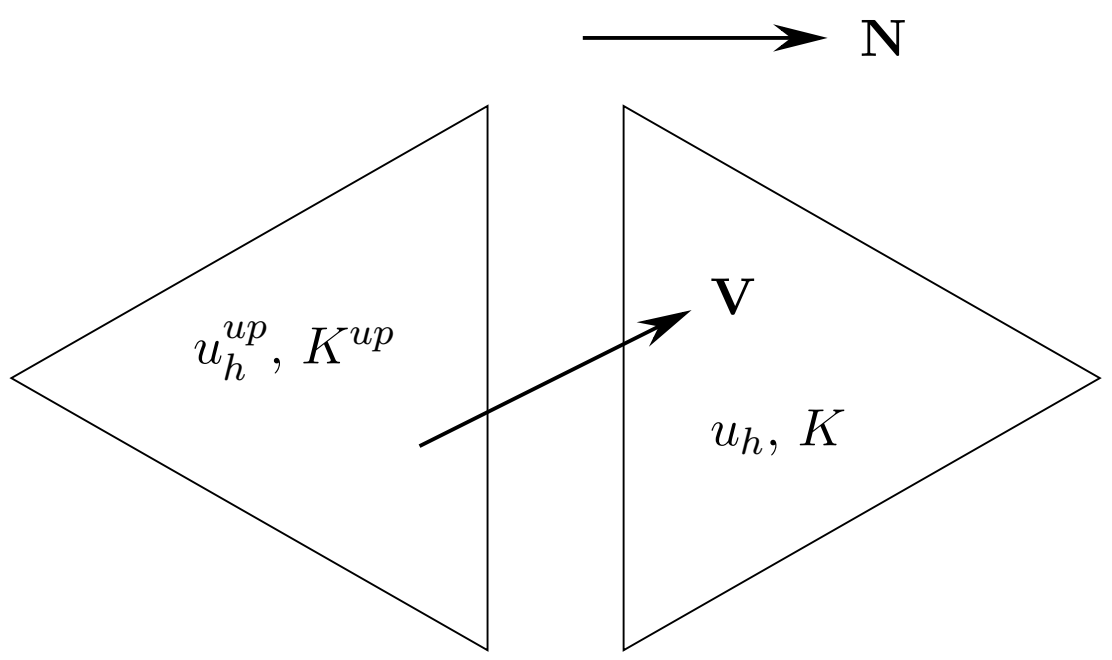

The proof will proceed in several steps.

{\em1. Coercivity:} The energy estimate follows directly upon
  substituting $(v_h, \vhat_h) = (u_h, \uhat_h)$ in the bilinear form.
  
{\em2. Monotonicity:} Write the bilinear form as
  \begin{multline*}
    a_{h\epsilon}
    \left((u_h,\uhat_h),(v_h,\vhat_h) \vph\right) 
    = \left(u_h,  (1-(1/2)\rmdiv(\bfV)) v_h \vph\right) 
    + \epsilon \langle u_h - \uhat_h, v_h - \vhat_h \rangle \\
    + \epsilon \left(\nabla\calL_h(\uhat_h) ,
      \nabla \calL_h(\vhat_h) \vph\right) 
    + \langle (u^{up}_h - u_h/2),  v_h (\bfV.\bfn) \rangle.
  \end{multline*}
  The matrix for the weighted inner product
  $\left(u_h, (1-(1/2)\rmdiv(\bfV)) v_h \vph\right)$ with piecewise constant
  functions is diagonal with positive diagonal terms, and since
  $$
  (u_h - \uhat_h)(v_h - \vhat_h)
  = (u_h, \uhat_h) \begin{bmatrix} 1 & -1 \\ -1 & 1 \end{bmatrix}
  \begin{pmatrix} v_h \\ \vhat_h \end{pmatrix},
  $$
  the penalty term can be written as a sum of M-matrices, so the
  sum of these two terms is monotone.

  The fact that the stiffness matrix of the bilinear form
  $(\nabla \calL_h(\uhat_h), \nabla \calL_h(\vhat_h))$ is an
  M-matrix when the mesh has acute simplices is well known and can be
  found in~\cite[Theorem~6.24]{barrenechea2025monotone}.

  The expression for the convective term can be written as
  \begin{multline*}
    \sum_{K \in \calT_h} \int_{\partial K} (u^{up}_h - u_h/2) v_h \bfV.\bfn
    = \sum_\Gammain \int_k (1/2) u_h v_h (-\bfV.\bfn) 
    + \sum_\Gammaout \int_k (1/2) u_h v_h (\bfV.\bfn) \\
    + \sum_{k \in \delTho} (1/2) \int_k \left(u^{up}_h v^{up}_h 
      - 2 u^{up}_h v_h + u_h v_h \vph\right) (\bfV.\bfN),
  \end{multline*}
  where the normal $\bfN$ on the internal interfaces is the one for
  which $\bfV.\bfN \geq 0$ and the downwind values are denoted $u_h$
  and $v_h$ (see Figure~\ref{fig:jumps}). Since $\bfV.\bfn < 0$ on
  $\Gammain$ and $\bfV.\bfn > 0$ on $\Gammaout$, the boundary terms
  are monotone. The interior terms are monotone by~\eqref{eqn:jumps-inequality} in Lemma~\ref{lem:interior-edges-jumps}.
  
{\em3. Consistency:} Fix $v \in \calD(\Omega)$ and let
  $(v_h, \vhat_h) = \iota_h(v)$ be the interpolant in Lemma~\ref{lem:interp}.
  From Lemma~\ref{lem:interp}, the following error estimates hold:
  \begin{align*}
    \ltwo{v-v_h} + \sqrt{h}  \langle v_h - \vhat_h\rangle
    & \leq C \ltwo{\nabla v} h, \\
    \ltwo{\nabla v - \nabla \calL_h(\vhat_h)}
    &\leq C \ltwo{D^2 v} h,
  \end{align*}
  so it suffices to verify the half-order rate for the convective term.
  The term arising from convection is bounded as
  $$
  \int_{\partial K} (v-v_h)^2 |\bfV.\bfn| 
  \leq \norm{v-v_h}^2_{L^4(\partial K)}
  \norm{\bfV}_{L^2(\partial K)} 
  \leq C \norm{\nabla v}^2_{L^4(K)} h^{3/2}_K
  \times \norm{\bfV}_{H^1(K)} / \sqrt{h_K}.
  $$
  Applying the Cauchy--Schwarz inequality to the sum over all elements
  bounds the right-hand side by $C \lfour{\nabla v}^2 \hone{\bfV} h$.
      
{\em4. Inflow Barrier:} Let $(\psi_h, \psihat_h) = \iota_h(\psi)$
  and fix $(v_h, \vhat_h) \in \calP_0(\calT_h) \times \calP_0(\delTh)$
  with $\vhat_h|_{\partial\Omega} = 0$. Since $v_h$ is piecewise
  constant, it is immediate that $(\psi_h, v_h) = (\psi,v_h)$. To
  bound the penalty term, Lemma~\ref{lem:interp} again gives
  $\langle \psi_h - \psihat_h\rangle \leq C \ltwo{\nabla \psi}
  \sqrt{h}$.
  Then
  $$
  \langle \psi_h - \psihat_h, v_h - \vhat_h \rangle
  \leq \langle \psi_h - \psihat_h \rangle \, \langle v_h - \vhat_h \rangle
  \leq C \langle v_h - \vhat_h \rangle \ltwo{\nabla \psi} \sqrt{h}.
  $$
  For the convective terms, integrating by parts and adding and
  subtracting $\psi$ gives
  \begin{eqnarray*}
    \int_{\partial K} (\psi^{up}_h - \psi_h) v_h \bfV.\bfn 
    &=& \int_K -\psi_h v_h \rmdiv(\bfV)
    + \int_{\partial K} \psi^{up}_h v_h \bfV.\bfn  \\
    &=& \int_K -\psi v_h \rmdiv(\bfV) + (\psi -\psi_h) v_h \rmdiv(\bfV)
    + \int_{\partial K} \psi^{up}_h v_h \bfV.\bfn  \\
    &=& \int_K \bfV.\nabla \psi \, v_h + (\psi -\psi_h) v_h \rmdiv(\bfV)
        + \int_{\partial K} (\psi^{up}_h-\psi) v_h \bfV.\bfn .
  \end{eqnarray*}
  Since the boundary integrand vanishes on $\partial \Omega$ and
  $\psi^{up}_h - \psi$ is single-valued on each edge/face, after summing over elements, we have
  \begin{multline*}
    \sum_{K \in \calT_h}
    \int_{\partial K} (\psi^{up}_h - \psi_h) v_h \bfV.\bfn 
    = \int_\Omega \bfV.\nabla \psi \, v_h  \\
       + \int_\Omega (\psi -\psi_h) v_h \rmdiv(\bfV)
       + \sum_{k \in \delTho}
        \int_k (\psi^{up}_h-\psi) [v_h] \bfV.\bfN . 
  \end{multline*}
  The first term on the second line is bounded by
  $C \ltwo{\nabla \psi} \ltwo{v_h} \linf{\rmdiv(\bfV)} h$, and the
    second term is bounded as
  \begin{eqnarray*}
    \sum_{k \in \delTho} \int_k (\psi^{up}_h-\psi) [v_h] \bfV.\bfN
    &\leq& \left(\sum_{k \in \delTho} \int_k [v_h]^2 |\bfV.\bfn| \right)^{1/2}
           \sum_{K \in \calT_h} \norm{\psi_h - \psi}_{L^4(\partial K)}
           \norm{\bfV}_{L^2(\partial K)}^{1/2} \\
    &\leq& C \norm{(v_h,\vhat_h)}_{h\epsilon} \lfour{\nabla \psi}
           \hone{\bfV}^{1/2} \sqrt{h}.        
  \end{eqnarray*}

{\em5. Outflow Cutoff:} Setting $(v_h, \vhat_h) = (\phi_h u_h,
  \phihat_h \uhat_h)$ with $(\phi_h, \phihat_h) = \iota_h(\phi)$ in
  the bilinear form gives
  \begin{multline*}
    a_{h\epsilon}\big((u_h,\uhat_h),
    (\phi_h u_h, \phihat_h \uhat_h) \big)
    = \int_\Omega \phi_h u_h^2
    + \epsilon \nabla \calL_h(\uhat_h).\nabla \calL_h(\phihat_h \uhat_h) \\
    + \sum_{K \in \calT_h} \int_{\partial K}
    \epsilon (u_h-\uhat_h)(\phi_h u_h - \phihat_h \uhat_h)
    + (u^{up}_h - u_h) \phi_h u_h \bfV.\bfn
    = \int_\Omega \phi_h f u_h.
  \end{multline*}
  Since $\phihat_h$ is the average of $\phi$ on each edge, the penalty
  term is evaluated as
  $$
  \langle u_h-\uhat_h, \phi_h u_h - \phihat_h \uhat_h \rangle
  = \langle u_h-\uhat_h, \phi (u_h - \uhat_h) \rangle
  + \langle u_h-\uhat_h, (\phi_h - \phi) u_h \rangle.
  $$
  The first term is $\langle \phi^{1/2} (u_h - \uhat_h) \rangle^2$
  and the second term is bounded as
  \begin{eqnarray*}
    \langle u_h-\uhat_h, (\phi_h - \phi) u_h \rangle
    &\leq& \langle u_h-\uhat_h \rangle
           \langle (\phi_h - \phi) u_h \rangle \\
    &\leq& C \langle u_h-\uhat_h \rangle \linf{\nabla \phi}
           \ltwo{u_h} \sqrt{h}.
  \end{eqnarray*}
  
  To identify the discrete product rule for the elliptic term, we
  repeatedly use an integration by parts identity for
  Crouzeix--Raviart functions together with the fact that
  $\hat{\phi}_h$ is the $L^2$-projection of $\phi$ onto
  $\calP_0(\delTh)$. For any $\hat{u}_h, \vhat_h\in \calP_0(\delTh)$
  and smooth $w$ with $(w_h, \hat{w}_h) = \iota_h(w)$, we have
 \begin{equation}\label{eqn:CR-int-by-parts}
  \int_K \nabla \calL_h(\vhat_h).\nabla \calL_h(\hat{w}_h \uhat_h)
    = \int_{\partial K} \hat{w}_h \uhat_h
        \nabla \calL_h(\vhat_h).\bfn  = \int_{\partial K} w \uhat_h
        \nabla \calL_h(\vhat_h).\bfn .
  \end{equation}
  We use~\eqref{eqn:CR-int-by-parts} with $\vhat_h=\uhat_h$
  and $w=\phi$, and add  $\pm (\phi-\phi_h)(\uhat_h-u_h)$ to
  get the following equality:
  \begin{align*}
    \int_K \nabla \calL_h(\uhat_h).\nabla \calL_h(\hat{\phi}_h \uhat_h)
    &= \int_{\partial K} \phi \uhat_h
      \nabla \calL_h(\uhat_h).\bfn \\
    &= \int_{\partial K}
      \left(\phi_h \uhat_h + u_h \phi - \phi_h u_h
      + (\phi - \phi_h)(\uhat_h - u_h) \right)
      \nabla \calL_h(\uhat_h).\bfn.
  \end{align*}
  Using the fact that $u_h$ and $\phi_h$ are constant on $K$ and applying
  the integration by parts identity~\eqref{eqn:CR-int-by-parts} to the
  first three terms on the right-hand side gives
  \begin{align*}
    \int_K \nabla \calL_h(\uhat_h).\nabla \calL_h(\hat{\phi}_h \uhat_h)  = 
    \int_K& \phi_h |\nabla \calL_h(\uhat_h)|^2
            + u_h \nabla \calL_h(\uhat_h) . \nabla \calL_h(\phi) 
    \\&+ \int_{\partial K}(\phi - \phi_h)(\uhat_h - u_h)
    \, \nabla \calL_h(\uhat_h).\bfn .
  \end{align*}
  Since $|\nabla \calL_h(\uhat_h)|^2$ is constant, $\phi$ may be
  substituted for $\phi_h$ in the first term. The last term
  (consistency error) is bounded using the Cauchy--Schwarz inequality, an error
  estimate for $\phi$, and a trace inequality:
  $$
  \sum_{K \in \calT_h} \int_{\partial K}(\phi - \phi_h)(\uhat_h - u_h)
  \, \nabla \calL_h(\uhat_h).\bfn
  \leq C \linf{\nabla \phi} \langle u_h - \uhat_h \rangle
  \ltwo{\nabla \calL_h(\uhat_h)} \sqrt{h}.
  $$
  To identify the convective term, add and subtract $\phi$:
  $$
  \sum_{K \in \calT_h}
  \int_{\partial K} (u^{up}_h - u_h) \phi_h u_h \bfV.\bfn
  = \sum_{K \in \calT_h} \int_{\partial K}
  (u^{up}_h - u_h) \phi u_h \bfV.\bfn
  + (u^{up}_h - u_h) (\phi_h - \phi) u_h \bfV.\bfn.
  $$
  The first term may be rewritten as
  \begin{eqnarray*}
    \sum_{K \in \calT_h} \int_{\partial K} (u^{up}_h - u_h) \phi u_h \bfV.\bfn
    &=& \sum_{K \in \calT_h} \int_{\partial K}
        (u^{up}_h - u_h/2) \phi u_h \bfV.\bfn
        - (u_h^2/2) \phi \bfV.\bfn \\
    &=& -\sum_{K \in \calT_h}\int_K \frac{1}{2}u_h^2 \, \rmdiv(\phi\bfV)
        + \int_{\partial K} (u^{up}_h - u_h/2) \phi u_h \bfV.\bfn \\
    &=& -\frac{1}{2}\int_\Omega u_h^2 \, \rmdiv(\phi\bfV)
        + \frac{1}{2}\int_{\Gammain} \phi u_h^2 |\bfV.\bfn|
       + \frac{1}{2}\sum_{k \in \delTho} \int_k \phi [u_h]^2 |\bfV.\bfn|.
  \end{eqnarray*}
  To derive the last line, note first that $\phi \bfV.\bfn$ vanishes
  on $\partial \Omega \setminus \Gammain$ and on $\Gammain$ it
  evaluates to $(\phi/2)u_h^2 |\bfV.\bfn|$ since $u^{up}_h \equiv 0$
  on the inflow boundary. Then, the internal faces/edges are handled
  using~\eqref{eqn:jumps-equality} in Lemma~\ref{lem:interior-edges-jumps}.

  To bound the consistency error of the convective term, bound the
  difference $\phi_h - \phi$ by $C \linf{\nabla \phi} h$ and the
  remaining terms with the Cauchy--Schwarz inequality to get
  \begin{eqnarray*}
    \lefteqn{ \sum_{K \in \calT_h} \int_{\partial K}
    (u^{up}_h - u_h) (\phi_h - \phi) u_h \bfV.\bfn } \\
    &\leq& 
    C \left(\sum_{K \in \calT_h} \int_{\partial K}
           (u^{up}_h - u_h)^2 |\bfV.\bfn| \right)^{1/2}
           \linf{\nabla \phi} h
           \left(\sum_{K \in \calT_h} \int_{\partial K}
           u_h^2 |\bfV.\bfn| \right)^{1/2} \\
    &\leq& 
    C \left(\ltwoin{u_h} + \ltwoDelTho{[u_h]} \right)
           \linf{\nabla \phi} h
           \left(\sum_{K \in \calT_h} 
           \norm{u_h}^2_{L^2(K)} \frac{1}{|K|}
           \int_{\partial K}
           |\bfV.\bfn| \right)^{1/2},
  \end{eqnarray*}
  where the last line follows since $u_h$ is constant on each element.
  Lemma~\ref{lem:trace} bounds each integral in the last line
  by $C h^{d/2}\sqrt{\log(h^{-1})} \hone{\bfV}$, and since $|K| \geq c h^d$, it follows that
  \begin{multline*}
    \sum_{K \in \calT_h} \int_{\partial K}
    (u^{up}_h - u_h) (\phi_h - \phi) u_h \bfV.\bfn \\
    \leq C \left(\ltwoin{u_h} + \ltwoDelTho{[u_h]} \right)
    \linf{\nabla \phi} \ltwo{u_h} \hone{\bfV}^{1/2} h^{1-d/4}(\log(h^{-1}))^{1/4}.
  \end{multline*}

\bibliographystyle{siam}
\bibliography{references}

\end{document}

%% file: commands.tex
\usepackage{amsmath,amssymb,amscd,subfigure}
\usepackage{amsfonts,graphicx,psfrag,esint}

\renewcommand{\Re}{\mathbb{R}}

\newcommand{\half}{\frac{1}{2}}
\newcommand{\weak}{\rightharpoonup}

\newcommand{\vph}{\vphantom{A^{A}_{A}}}
\newcommand{\sst}{\,\mid\,}

\newcommand{\dbydp}[2]{\frac{\partial #1}{\partial #2}}

\newcommand{\uhat}{\hat{u}}
\newcommand{\vhat}{\hat{v}}

\newcommand{\Phat}{\hat{P}}

\newcommand{\phihat}{\hat{\phi}}

\newcommand{\psihat}{\hat{\psi}}

\newcommand{\Omegabar}{\bar{\Omega}}

\newcommand{\ut}{\tilde{u}}

\newcommand{\bfe}{{\bf e}}

\newcommand{\bfn}{{\bf n}}

\newcommand{\bfx}{{\bf x}}

\newcommand{\bfN}{{\bf N}}

\newcommand{\bfV}{{\bf V}}

\newcommand{\calD}{{\cal D}}

\newcommand{\calL}{{\cal L}}

\newcommand{\calP}{{\cal P}}

\newcommand{\calT}{{\cal T}}
\newcommand{\calU}{{\cal U}}

\newcommand{\norm}[1]{\| {#1} \|}
\newcommand{\Hone}{{H^1(\Omega)}}
\newcommand{\hone}[1]{\norm{#1}_\Hone}
\newcommand{\Honeo}{{H^1_0(\Omega)}}

\newcommand{\Htwo}{{H^2(\Omega)}}
\newcommand{\htwo}[1]{\norm{#1}_\Htwo}

\newcommand{\Lone}{{L^1(\Omega)}}

\newcommand{\Ltwo}{{L^2(\Omega)}}
\newcommand{\ltwo}[1]{\norm{#1}_\Ltwo}

\newcommand{\Linf}{{L^\infty(\Omega)}}
\newcommand{\linf}[1]{\norm{#1}_{\Linf}}
\newcommand{\Lp}{{L^p(\Omega)}}
\newcommand{\lp}[1]{\norm{#1}_\Lp}

\newcommand{\Wonep}{{W^{1,p}(\Omega)}}

\newcommand{\Wonepo}{{W^{1,p}_0(\Omega)}}

\newcommand{\Woneinf}{{W^{1,\infty}(\Omega)}}

\newcommand{\Lq}{{L^q(\Omega)}}

\newcommand{\Lfour}{{L^4(\Omega)}}
\newcommand{\lfour}[1]{\norm{#1}_\Lfour}
\newcommand{\Wonefour}{{W^{1,4}(\Omega)}}

%% file: jumps.tex
\begin{tikzpicture}[scale=2]
  % two triangles and that share an edge
  
  \draw (-2.0,0.0) -- (-0.25, 1.0) -- (-0.25,-1.0) -- cycle;
  \draw ( 2.0,0.0) -- ( 0.25, 1.0) -- ( 0.25,-1.0) -- cycle;

  \node[right] at (-1.4,  0.1) {$u^{up}_h$, $K^{up}$};
  \node[right] at ( 0.5, -0.2) {$u_h$, $K$};

  % \draw[thick] (0.0,-1.0) -- (0.0, 1.0);  % edge

  % velocity arrow
  \draw[thick,-{Stealth[length=4mm, width=2mm]}] (-0.5, -0.25) -- (0.5,0.25);
  \node[right] at (0.5, 0.30) {$\bfV$};

  % normal arrow 
  \draw[thick,-{Stealth[length=4mm, width=2mm]}] (0.1, 1.25) -- (1.0,1.25);
  \node[right] at (1.05,1.25) {$\bfN$};

\end{tikzpicture}

%% file: paperII.bbl
\begin{thebibliography}{10}

\bibitem{aguillon2018error}
{\sc N.~Aguillon and F.~Boyer}, {\em Error estimate for the upwind scheme for
  the linear transport equation with boundary data}, IMA Journal of Numerical
  Analysis, 38 (2018), pp.~669--719.

\bibitem{armitage1985convexity}
{\sc D.~H. Armitage and {\"U}.~Kuran}, {\em The convexity of a domain and the
  superharmonicity of the signed distance function}, Proceedings of the
  American Mathematical Society, 93 (1985), pp.~598--600.

\bibitem{ayuso2009discontinuous}
{\sc B.~Ayuso and L.~D. Marini}, {\em Discontinuous {G}alerkin methods for
  advection-diffusion-reaction problems}, SIAM Journal on Numerical Analysis,
  47 (2009), pp.~1391--1420.

\bibitem{babuvska1992locking}
{\sc I.~Babu{\v{s}}ka and M.~Suri}, {\em On locking and robustness in the
  finite element method}, SIAM Journal on Numerical Analysis, 29 (1992),
  pp.~1261--1293.

\bibitem{bacuta2026explaining}
{\sc C.~Bacuta}, {\em Explaining oscillatory behavior in convection--diffusion
  discretization}, 2026.
\newblock arXiv:2601.09657.

\bibitem{barrenechea2025monotone}
{\sc G.~R. Barrenechea, V.~John, and P.~Knobloch}, {\em Monotone
  discretizations for elliptic second order partial differential equations},
  vol.~61 of Springer Series in Computational Mathematics, Springer, Cham,
  2025.

\bibitem{bause2002uniform}
{\sc M.~Bause and P.~Knabner}, {\em Uniform error analysis for
  {L}agrange--{G}alerkin approximations of convection-dominated problems}, SIAM
  Journal on Numerical Analysis, 39 (2002), pp.~1954--1984.

\bibitem{boyer2005trace}
{\sc F.~Boyer}, {\em Trace theorems and spatial continuity properties for the
  solutions of the transport equation}, Differential and Integral Equations, 18
  (2005), pp.~891--934.

\bibitem{boyer2012analysis}
\leavevmode\vrule height 2pt depth -1.6pt width 23pt, {\em Analysis of the
  upwind finite volume method for general initial- and boundary-value transport
  problems}, IMA Journal of Numerical Analysis, 32 (2012), pp.~1404--1439.

\bibitem{boyer2012mathematical}
{\sc F.~Boyer and P.~Fabrie}, {\em Mathematical Tools for the Study of the
  Incompressible Navier-Stokes Equations and Related Models}, vol.~183 of
  Applied Mathematical Sciences, Springer, New York, 2013.

\bibitem{brezis1972convergence}
{\sc H.~Br{\'e}zis and A.~Pazy}, {\em Convergence and approximation of
  semigroups of nonlinear operators in banach spaces}, Journal of Functional
  Analysis, 9 (1972), pp.~63--74.

\bibitem{cockburn2008optimal}
{\sc B.~Cockburn, B.~Dong, and J.~Guzm{\'a}n}, {\em Optimal convergence of the
  original {DG} method for the transport-reaction equation on special meshes},
  SIAM Journal on Numerical Analysis, 46 (2008), pp.~1250--1265.

\bibitem{crandall1971generation}
{\sc M.~G. Crandall and T.~M. Liggett}, {\em Generation of semi-groups of
  nonlinear transformations on general {B}anach spaces}, American Journal of
  Mathematics, 93 (1971), pp.~265--298.

\bibitem{devinatz1974asymptotic}
{\sc A.~Devinatz, R.~Ellis, and A.~Friedman}, {\em The asymptotic behavior of
  the first real eigenvalue of second order elliptic operators with a small
  parameter in the highest derivatives, {II}}, Indiana University Mathematics
  Journal, 23 (1974), pp.~991--1011.

\bibitem{di2011mathematical}
{\sc D.~A. Di~Pietro and A.~Ern}, {\em Mathematical aspects of discontinuous
  Galerkin methods}, vol.~69 of Math{\'e}matiques \& Applications, Springer,
  Berlin, Heidelberg, 2012.

\bibitem{diperna1989ordinary}
{\sc R.~J. DiPerna and P.-L. Lions}, {\em Ordinary differential equations,
  transport theory and {S}obolev spaces}, Inventiones Mathematicae, 98 (1989),
  pp.~511--547.

\bibitem{egger2010hybrid}
{\sc H.~Egger and J.~Sch{\"o}berl}, {\em A hybrid mixed discontinuous
  {G}alerkin finite-element method for convection--diffusion problems}, IMA
  Journal of Numerical Analysis, 30 (2010), pp.~1206--1234.

\bibitem{ern2006discontinuous}
{\sc A.~Ern and J.-L. Guermond}, {\em Discontinuous {G}alerkin methods for
  {F}riedrichs' systems. {I}. {G}eneral theory}, SIAM Journal on Numerical
  Analysis, 44 (2006), pp.~753--778.

\bibitem{fu2015analysis}
{\sc G.~Fu, W.~Qiu, and W.~Zhang}, {\em An analysis of {HDG} methods for
  convection-dominated diffusion problems}, ESAIM: Mathematical Modelling and
  Numerical Analysis, 49 (2015), pp.~225--256.

\bibitem{itokappel1998trotter}
{\sc K.~Ito and F.~Kappel}, {\em The {Trotter--Kato} theorem and approximation
  of {PDEs}}, Mathematics of Computation, 67 (1998), pp.~21--44.

\bibitem{jiang2025stability}
{\sc J.~Jiang, N.~J. Walkington, and Y.~Yue}, {\em Stability and convergence of
  {HDG} schemes under minimal regularity}, SIAM Journal on Numerical Analysis,
  63 (2025), pp.~2048--2071.

\bibitem{jin1999efficient}
{\sc S.~Jin}, {\em Efficient asymptotic-preserving ({AP}) schemes for some
  multiscale kinetic equations}, SIAM Journal on Scientific Computing, 21
  (1999), pp.~441--454.

\bibitem{johnson1986analysis}
{\sc C.~Johnson and J.~Pitk{\"a}ranta}, {\em An analysis of the discontinuous
  {G}alerkin method for a scalar hyperbolic equation}, Mathematics of
  Computation, 46 (1986), pp.~1--26.

\bibitem{peterson1991note}
{\sc T.~E. Peterson}, {\em A note on the convergence of the discontinuous
  {G}alerkin method for a scalar hyperbolic equation}, SIAM Journal on
  Numerical Analysis, 28 (1991), pp.~133--140.

\bibitem{riviere2008discontinuous}
{\sc B.~Rivi{\`e}re}, {\em Discontinuous Galerkin methods for solving elliptic
  and parabolic equations: theory and implementation}, SIAM, 2008.

\bibitem{roos2008robust}
{\sc H.-G. Roos, M.~Stynes, and L.~Tobiska}, {\em Robust numerical methods for
  singularly perturbed differential equations: convection-diffusion-reaction
  and flow problems}, vol.~24 of Springer Series in Computational Mathematics,
  Springer, Berlin, second~ed., 2008.

\bibitem{schlichting2018analysis}
{\sc A.~Schlichting and C.~Seis}, {\em Analysis of the implicit upwind finite
  volume scheme with rough coefficients}, Numerische Mathematik, 139 (2018),
  pp.~155--186.

\bibitem{schoberl2014c++}
{\sc J.~Sch{\"o}berl}, {\em C++ 11 implementation of finite elements in
  ngsolve}, Institute for analysis and scientific computing, Vienna University
  of Technology, 30 (2014).

\bibitem{tian2014asymptotically}
{\sc X.~Tian and Q.~Du}, {\em Asymptotically compatible schemes and
  applications to robust discretization of nonlocal models}, SIAM Journal on
  Numerical Analysis, 52 (2014), pp.~1641--1665.

\bibitem{walkington2005convergence}
{\sc N.~J. Walkington}, {\em Convergence of the discontinuous {G}alerkin method
  for discontinuous solutions}, SIAM Journal on Numerical Analysis, 42 (2005),
  pp.~1801--1817.

\end{thebibliography}
